\documentclass[12pt]{amsart}
\usepackage{amsmath}
\usepackage{amssymb,latexsym}
\usepackage{amscd}
\usepackage{comment}
\usepackage{xspace}
\usepackage{verbatim}
\usepackage[english]{babel}

\begin{document}

\numberwithin{equation}{section}
\title[Elliptic equations with nonlinear absorption]{Elliptic equations with nonlinear absorption depending on the solution and its gradient}
\author{Moshe Marcus }
\address{Department of Mathematics, Technion\\
 Haifa 32000, ISRAEL}
 \email{marcusm@math.technion.ac.il}
\author{Phuoc-Tai Nguyen}
\address{Department of Mathematics, Technion\\
 Haifa 32000, ISRAEL}
 \email{nguyenphuoctai.hcmup@gmail.com}

\date{}

 \maketitle


\newcommand{\txt}[1]{\;\text{ #1 }\;}
\newcommand{\tbf}{\textbf}
\newcommand{\tit}{\textit}
\newcommand{\tsc}{\textsc}
\newcommand{\trm}{\textrm}
\newcommand{\mbf}{\mathbf}
\newcommand{\mrm}{\mathrm}
\newcommand{\bsym}{\boldsymbol}
\newcommand{\scs}{\scriptstyle}
\newcommand{\sss}{\scriptscriptstyle}
\newcommand{\txts}{\textstyle}
\newcommand{\dsps}{\displaystyle}
\newcommand{\fnz}{\footnotesize}
\newcommand{\scz}{\scriptsize}
\newcommand{\be}{\begin{equation}}
\newcommand{\bel}[1]{\begin{equation}\label{#1}}
\newcommand{\ee}{\end{equation}}
\newcommand{\eqnl}[2]{\begin{equation}\label{#1}{#2}\end{equation}}
\newcommand{\barr}{\begin{eqnarray}}
\newcommand{\earr}{\end{eqnarray}}
\newcommand{\bars}{\begin{eqnarray*}}
\newcommand{\ears}{\end{eqnarray*}}
\newcommand{\nnu}{\nonumber \\}
\newtheorem{subn}{\name}
\renewcommand{\thesubn}{}
\newcommand{\bsn}[1]{\def\name{#1}\begin{subn}}
\newcommand{\esn}{\end{subn}}
\newtheorem{sub}{\name}[section]
\newcommand{\dn}[1]{\def\name{#1}}   
\newcommand{\bs}{\begin{sub}}
\newcommand{\es}{\end{sub}}
\newcommand{\bsl}[1]{\begin{sub}\label{#1}}
\newcommand{\bth}[1]{\def\name{Theorem}
\begin{sub}\label{t:#1}}
\newcommand{\blemma}[1]{\def\name{Lemma}
\begin{sub}\label{l:#1}}
\newcommand{\bcor}[1]{\def\name{Corollary}
\begin{sub}\label{c:#1}}
\newcommand{\bdef}[1]{\def\name{Definition}
\begin{sub}\label{d:#1}}
\newcommand{\bprop}[1]{\def\name{Proposition}
\begin{sub}\label{p:#1}}
\newcommand{\R}{\eqref}
\newcommand{\rth}[1]{Theorem~\ref{t:#1}}
\newcommand{\rlemma}[1]{Lemma~\ref{l:#1}}
\newcommand{\rcor}[1]{Corollary~\ref{c:#1}}
\newcommand{\rdef}[1]{Definition~\ref{d:#1}}
\newcommand{\rprop}[1]{Proposition~\ref{p:#1}}
\newcommand{\BA}{\begin{array}}
\newcommand{\EA}{\end{array}}
\newcommand{\BAN}{\renewcommand{\arraystretch}{1.2}
\setlength{\arraycolsep}{2pt}\begin{array}}
\newcommand{\BAV}[2]{\renewcommand{\arraystretch}{#1}
\setlength{\arraycolsep}{#2}\begin{array}}
\newcommand{\BSA}{\begin{subarray}}
\newcommand{\ESA}{\end{subarray}}
\newcommand{\BAL}{\begin{aligned}}
\newcommand{\EAL}{\end{aligned}}
\newcommand{\BALG}{\begin{alignat}}
\newcommand{\EALG}{\end{alignat}}
\newcommand{\BALGN}{\begin{alignat*}}
\newcommand{\EALGN}{\end{alignat*}}
\newcommand{\note}[1]{\textit{#1.}\hspace{2mm}}
\newcommand{\Proof}{\noindent \note{\it Proof}}
\newcommand{\Remark}{\noindent \note{\it Remark}}
\newcommand{\modin}{$\,$\\[-4mm] \indent}
\newcommand{\forevery}{\quad \forall}
\newcommand{\set}[1]{\{#1\}}
\newcommand{\setdef}[2]{\{\,#1:\,#2\,\}}
\newcommand{\setm}[2]{\{\,#1\mid #2\,\}}
\newcommand{\mt}{\mapsto}
\newcommand{\lra}{\longrightarrow}
\newcommand{\lla}{\longleftarrow}
\newcommand{\llra}{\longleftrightarrow}
\newcommand{\Lra}{\Longrightarrow}
\newcommand{\Lla}{\Longleftarrow}
\newcommand{\Llra}{\Longleftrightarrow}
\newcommand{\warrow}{\rightharpoonup}
\newcommand{
\paran}[1]{\left (#1 \right )}
\newcommand{\sqbr}[1]{\left [#1 \right ]}
\newcommand{\curlybr}[1]{\left \{#1 \right \}}
\newcommand{\abs}[1]{\left |#1\right |}
\newcommand{\norm}[1]{\left \|#1\right \|}
\newcommand{
\paranb}[1]{\big (#1 \big )}
\newcommand{\lsqbrb}[1]{\big [#1 \big ]}
\newcommand{\lcurlybrb}[1]{\big \{#1 \big \}}
\newcommand{\absb}[1]{\big |#1\big |}
\newcommand{\normb}[1]{\big \|#1\big \|}
\newcommand{
\paranB}[1]{\Big (#1 \Big )}
\newcommand{\absB}[1]{\Big |#1\Big |}
\newcommand{\normB}[1]{\Big \|#1\Big \|}
\newcommand{\produal}[1]{\langle #1 \rangle}

\newcommand{\thkl}{\rule[-.5mm]{.3mm}{3mm}}
\newcommand{\thknorm}[1]{\thkl #1 \thkl\,}
\newcommand{\trinorm}[1]{|\!|\!| #1 |\!|\!|\,}
\newcommand{\bang}[1]{\langle #1 \rangle}
\def\angb<#1>{\langle #1 \rangle}
\newcommand{\vstrut}[1]{\rule{0mm}{#1}}
\newcommand{\rec}[1]{\frac{1}{#1}}
\newcommand{\opname}[1]{\mbox{\rm #1}\,}
\newcommand{\supp}{\opname{supp}}
\newcommand{\dist}{\opname{dist}}
\newcommand{\myfrac}[2]{{\displaystyle \frac{#1}{#2} }}
\newcommand{\myint}[2]{{\displaystyle \int_{#1}^{#2}}}
\newcommand{\mysum}[2]{{\displaystyle \sum_{#1}^{#2}}}
\newcommand {\dint}{{\displaystyle \myint\!\!\myint}}
\newcommand{\q}{\quad}
\newcommand{\qq}{\qquad}
\newcommand{\hsp}[1]{\hspace{#1mm}}
\newcommand{\vsp}[1]{\vspace{#1mm}}
\newcommand{\ity}{\infty}
\newcommand{\prt}{\partial}
\newcommand{\sms}{\setminus}
\newcommand{\ems}{\emptyset}
\newcommand{\ti}{\times}
\newcommand{\pr}{^\prime}
\newcommand{\ppr}{^{\prime\prime}}
\newcommand{\tl}{\tilde}
\newcommand{\sbs}{\subset}
\newcommand{\sbeq}{\subseteq}
\newcommand{\nind}{\noindent}
\newcommand{\ind}{\indent}
\newcommand{\ovl}{\overline}
\newcommand{\unl}{\underline}
\newcommand{\nin}{\not\in}
\newcommand{\pfrac}[2]{\genfrac{(}{)}{}{}{#1}{#2}}

\def\ga{\alpha}     \def\gb{\beta}       \def\gg{\gamma}
\def\gc{\chi}       \def\gd{\delta}      \def\ge{\epsilon}
\def\gth{\theta}                         \def\vge{\varepsilon}
\def\gf{\phi}       \def\vgf{\varphi}    \def\gh{\eta}
\def\gi{\iota}      \def\gk{\kappa}      \def\gl{\lambda}
\def\gm{\mu}        \def\gn{\nu}         \def\gp{\pi}
\def\vgp{\varpi}    \def\gr{\rho}        \def\vgr{\varrho}
\def\gs{\sigma}     \def\vgs{\varsigma}  \def\gt{\tau}
\def\gu{\upsilon}   \def\gv{\vartheta}   \def\gw{\omega}
\def\gx{\xi}        \def\gy{\psi}        \def\gz{\zeta}
\def\Gg{\Gamma}     \def\Gd{\Delta}      \def\Gf{\Phi}
\def\Gth{\Theta}
\def\Gl{\Lambda}    \def\Gs{\Sigma}      \def\Gp{\Pi}
\def\Gw{\Omega}     \def\Gx{\Xi}         \def\Gy{\Psi}

\def\CS{{\mathcal S}}   \def\CM{{\mathcal M}}   \def\CN{{\mathcal N}}
\def\CR{{\mathcal R}}   \def\CO{{\mathcal O}}   \def\CP{{\mathcal P}}
\def\CA{{\mathcal A}}   \def\CB{{\mathcal B}}   \def\CC{{\mathcal C}}
\def\CD{{\mathcal D}}   \def\CE{{\mathcal E}}   \def\CF{{\mathcal F}}
\def\CG{{\mathcal G}}   \def\CH{{\mathcal H}}   \def\CI{{\mathcal I}}
\def\CJ{{\mathcal J}}   \def\CK{{\mathcal K}}   \def\CL{{\mathcal L}}
\def\CT{{\mathcal T}}   \def\CU{{\mathcal U}}   \def\CV{{\mathcal V}}
\def\CZ{{\mathcal Z}}   \def\CX{{\mathcal X}}   \def\CY{{\mathcal Y}}
\def\CW{{\mathcal W}} \def\CQ{{\mathcal Q}}
\def\BBA {\mathbb A}   \def\BBb {\mathbb B}    \def\BBC {\mathbb C}
\def\BBD {\mathbb D}   \def\BBE {\mathbb E}    \def\BBF {\mathbb F}
\def\BBG {\mathbb G}   \def\BBH {\mathbb H}    \def\BBI {\mathbb I}
\def\BBJ {\mathbb J}   \def\BBK {\mathbb K}    \def\BBL {\mathbb L}
\def\BBM {\mathbb M}   \def\BBN {\mathbb N}    \def\BBO {\mathbb O}
\def\BBP {\mathbb P}   \def\BBR {\mathbb R}    \def\BBS {\mathbb S}
\def\BBT {\mathbb T}   \def\BBU {\mathbb U}    \def\BBV {\mathbb V}
\def\BBW {\mathbb W}   \def\BBX {\mathbb X}    \def\BBY {\mathbb Y}
\def\BBZ {\mathbb Z}

\def\GTA {\mathfrak A}   \def\GTB {\mathfrak B}    \def\GTC {\mathfrak C}
\def\GTD {\mathfrak D}   \def\GTE {\mathfrak E}    \def\GTF {\mathfrak F}
\def\GTG {\mathfrak G}   \def\GTH {\mathfrak H}    \def\GTI {\mathfrak I}
\def\GTJ {\mathfrak J}   \def\GTK {\mathfrak K}    \def\GTL {\mathfrak L}
\def\GTM {\mathfrak M}   \def\GTN {\mathfrak N}    \def\GTO {\mathfrak O}
\def\GTP {\mathfrak P}   \def\GTR {\mathfrak R}    \def\GTS {\mathfrak S}
\def\GTT {\mathfrak T}   \def\GTU {\mathfrak U}    \def\GTV {\mathfrak V}
\def\GTW {\mathfrak W}   \def\GTX {\mathfrak X}    \def\GTY {\mathfrak Y}
\def\GTZ {\mathfrak Z}   \def\GTQ {\mathfrak Q}

\def\sign{\mathrm{sign\,}}
\def\bdw{\prt\Gw}
\def\nabu{|\nabla u|}
\def\tr{\mathrm{tr\,}}
\def\req{\eqref}
\def\sth{such that\xspace}
\tableofcontents

\begin{abstract}
We study positive solutions of equation (E1) $-\Delta u + u^p|\nabla u|^q= 0$ ($0\leq p$, $0\leq q\leq 2$, $p+q>1$) and (E2) $-\Delta u + u^p + |\nabla u|^q =0$ ($p>1$, $1<q\leq 2$)  in a smooth bounded domain $\Omega \subset \mathbb{R}^N$. We obtain a sharp condition on $p$ and $q$ under which, for every positive, finite Borel measure $\mu$ on $\partial \Omega$,  there exists a solution such that $u=\mu$ on $\partial \Omega$. Furthermore, if the condition mentioned above fails then any isolated point singularity on $\partial \Omega$ is removable, namely there is no positive solution that vanishes on $\partial \Omega$ everywhere except at one point. With respect to (E2) we also prove uniqueness and discuss solutions that blow-up on a compact subset of $\partial \Omega$. In both cases we obtain a classification of positive solutions with an isolated boundary singularity. Finally, in Appendix A a uniqueness result for a class of quasilinear equations is provided. This class includes (E1) when $p=0$ but not the general case.
\end{abstract} \medskip

 \smallskip


\medskip

\section{Introduction}
In this paper, we are concerned with the boundary value problems with measure data for equations of the form
\bel{A0} -\Gd u + H(x,u,\nabla u) = 0 \ee
in $\Gw$ where $\Gw$ is a $C^2$ bounded domain in $\BBR^N$ and $H \in C(\Gw \times \BBR \times \BBR^N)$, $H\geq 0$.

The case where $H$ depends only on $u$, has been intensively studied, especially the following typical equation
\bel{Ap} -\Gd u + |u|^p\sign u = 0  \ee
with $p>1$ (see Dynkin \cite{Dybook1, Dybook2}, Le Gall \cite{Lg2}, Gmira and V\'eron \cite{GV}, Marcus and V\'eron \cite{MV1, MV3, MVbook}, Marcus\cite{Ma} and the references therein). In \cite{GV} it was shown that  \eqref{Ap} admits a {\it critical value}
\bel{pc} p_c=\frac{N+1}{N-1}. \ee
such that, for $1<p<p_c$, the boundary value problem
\begin{equation}\label{u^p} \left\{ \BAL
 -\Gd u+|u|^p\sign u &= 0 \q \text{in}\;\Gw\\
 u&=\mu \q \text{on}\;\prt\Gw
\EAL \right. \end{equation}
  has a unique solution for every $\mu\in \GTM(\bdw)$ (= space of finite Borel measures on $\prt\Gw$). The boundary data is attained as a weak limit of measures.
 Moreover isolated boundary singularities of solutions of \eqref{Ap} can be completely described. For more general results on positive solutions of (1.2) with singular sets on the boundary see \cite{MV3,MV4}. For a treatment of more general equations (where the absorption term $H$ depends on $(x,u)$) see \cite{AnMa}.

The case where $H$ depends only on $\nabla u$ has been recently investigated by P.T. Nguyen and L. V\'eron \cite{NV}. For equations of the form
	\bel{Aq} -\Gd u + g(|\nabla u|) = 0 \q\text{in}\;\Gw \ee
they obtained a sufficient conditions on $g$ in order that the boundary value problem for \eqref{Aq} with measure boundary data  have  a solution for every measure in $\GTM(\bdw)$.
If the nonlinearity is of power type, namely $g(|\nabla u|)=|\nabla u|^q$ with $1 \leq q \leq 2$, they showed that the critical value for \eqref{Aq} is
\bel{qc} q_c=\frac{N+1}{N} \ee
and, for $1<q<q_c$, they provided a complete description of the positive solutions with isolated singularities on the boundary. The question of uniqueness for \eqref{Aq} and some related equations in subcritical case is treated in Appendix A of the present paper by the second author. The proof is based on  a technique of \cite{Po} adapted to the present case.

Notice that when $q >2$, by \cite{Li1}  if $u \in C^2(\Gw)$ is a positive solution of \eqref{Aq} then $u$ is bounded in $\Gw$. Therefore solutions may exist only for boundary data represented by a bounded function.

In the present paper, we study boundary value problems and boundary singularities of positive solutions of \eqref{A0} when {\it $H$ depends on both $u$ and $\nabla u$}. It is convenient to use the following notation: $H\circ u$ is the function given by
 $$(H\circ u)(x)=H(x,u(x),\!\nabla u(x)).$$

 We study the case of subquadratic growth in the gradient and concentrate on two model cases:
\bel{multi} H(x,t,\gx) =   t^p|\gx|^q \ee
where $p>0$, $0 \leq q \leq 2$ and
\bel{add}  H(x,t,\gx)  = t^p + |\gx|^q  \ee
where $p \geq 1$, $1 \leq q \leq 2$.

Equation \eqref{A0} with $H$ as in \eqref{add} was studied in \cite{AGMQ}, \cite{BaGi}; existence and uniqueness of large solutions was established  when $1<p<q\leq 2$. When $H$ is given by \eqref{multi}, there exists no large solution of equation \eqref{A0}. To our knowledge, up to now, there is no publication treating boundary value problems with measure data for \eqref{A0} and $H$ as in \eqref{multi} or \eqref{add}.

 The main difficulty that one encounters in the study of these problems: the inequality $u \leq v$ does not imply any relation between $|\nabla u|$ and $|\nabla v|$. Moreover, in general, the sum of two supersolutions of \eqref{A0} is not a supersolution. In addition, for $H$ as in  \eqref{multi} there is no a priori estimate  of solutions of \eqref{A0} or of their gradient. (However an  upper estimate is available for families of solutions satisfying certain auxiliary conditions.) On the other hand, when $H$ satisfies \eqref{multi} equation \eqref{A0} admits a similarity transformation; when $H$ is as in \eqref{add}, the equation does not admit a similarity transformation unless $p = \frac{q}{2-q}$.
\medskip

Before stating our main results we introduce some definitions.
\bdef{sol} (i) A function $u$ is  a (weak) solution of \eqref{A0} if $u \in L_{loc}^1(\Gw)$, $H \circ u \in L_{loc}^1(\Gw)$ and $u$ satisfies \eqref{A0} in the sense of distribution.

(ii)  Let $\mu\in \GTM(\prt\Gw)$. A function $u$ is  a solution of
\bel{A1} \left\{ \BAL -\Gd u +  H\circ u&= 0 \q  \text{in } \Gw \\
u&=\mu \q  \text{on } \prt\Gw
\EAL \right. \ee
if $u$ satisfies the equation and has boundary trace $\mu$ (see \rdef{Mtrace}).
\es

\noindent\note{Remark} It can be shown  that (see \rth{moderate} below) $u$ is a solution of \eqref{A1} if and only if $u \in L^1(\Gw)$, $H\circ u \in L^1_\rho(\Gw)$ and $u$ satisfies
\bel{A2}
\int_{\Gw}\left(-u\Gd\gz+(H\circ u)\gz\right)dx=-\int_{\prt\Gw}\frac{\prt \gz}{\prt \bf n}d\mu \forevery \gz \in C_0^2(\ovl \Gw)
\ee
where $\rho(x)=\dist(x,\prt \Gw)$, $\bf n$ denotes the  outward normal unit vector on $\prt \Gw$ and $C_0^2(\ovl \Gw)=\{u\in C^2(\ovl \Gw): u=0\q \text{on }\prt\Gw\}$.

\bdef{moderate} A positive solution of \eqref{A0} is moderate if $H\circ u\in L^1_\gr(\Gw)$.
\es

\bdef{subcriticality} A nonlinearity $H$ is called \emph{subcritical} if the problem \eqref{A1} admits a solution for every positive bounded measure $\gm$ on $\prt \Gw$. Otherwise, $H$ is called \emph{supercritical}. \es

Put
\bel{mpq} m_{p,q}:=\max\left\{p,\frac{q}{2-q}\right\}. \ee
The first theorem provides a sufficient condition for $H$ to be subcritical and a stability result relative to weak convergence of data. As shown later on (see Theorem F) the sufficient condition is also necessary for
subcriticality of $H$.
 \medskip

\noindent{\bf Theorem A.} {\it Assume either $H$ satisfies \eqref{multi} with $0<N(p+q-1)<p+1$ or  \eqref{add} with $m_{p,q}<p_c$. Then $H$ is subcritical and the following stability result holds:

Let $\{\gm_n\}$ be a sequence of positive finite measures on $\prt\Gw$ converging weakly to a positive finite measure $\gm$  and $\{u_{\gm_n}\}$ be a sequence of corresponding solutions of \eqref {A1} with $\gm=\gm_n$. Then there exists a subsequence such that $\{ u_{\gm_{n_k}}\}$ converges to a solution $ u_\gm$ of \eqref {A1} in $L^1(\Gw)$ and $\{H\circ u_{\gm_{n_k}}\}$ converges to $H \circ u$ in $L^1_\rho(\Gw)$.} \medskip

\Remark The method of proof of this theorem is classical. It is based on estimates in weak $L^p$ space and compactness of approximating solutions. The results stated in Theorem A can be extended, in the same way, to the following cases:
\bel{multia} 0 \leq H(x,t,\gx) \leq a_1(x)t^p|\gx|^q  \forevery (x,t,\gx) \in \Gw \ti \BBR_+ \ti \BBR^N \ee
where $p \geq 0$, $q \geq 0$,  $0<N(p+q-1)<p+1$ , $a_1 \in L^\infty(\Gw)$ and $a_1>c>0$;
\bel{adda} 0 \leq H(x,t,\gx)  \leq a_2(x)f(t) + a_3(x)g(|\gx|)  \forevery (x,t,\gx) \in \Gw \ti \BBR_+ \ti \BBR^N \ee
where $a_i \in L^\infty(\Gw)$, $a_i>c>0$ ($i=2,3$), $f$ and $g$ are positive, nondecreasing, continuous functions in $\BBR_+$, satisfying $f(0)=g(0)=0$ and
$$ \myint{1}{\infty}t^{-\frac{2N}{N-1}}f(t)dt<\infty, \q \myint{1}{\infty}t^{-\frac{2N+1}{N}}g(t)dt<\infty. $$

The next theorem presents an  uniqueness result when $H$ satisfies \eqref{add}. \medskip

\noindent{\bf Theorem B.} {\it Assume that $H$ satisfies \eqref{add} and $m_{p,q}<p_c$. Then \eqref{A1} has a \emph{unique} solution for every $\gm \in \GTM^+(\prt \Gw)$.} \medskip

The uniqueness of solutions of problem \eqref{A1}  when $H$ satisfies \eqref{multi} ($0<N(p+q-1)<p+1$) remains open. However we establish uniqueness in the case that $\gm$ is concentrated at a point.

In the next theorems we discuss solutions with an isolated singularity at a point $A\in  \bdw$. Without loss of generality we assume that $A$ is the origin. \medskip

\noindent{\bf Theorem C.} {\it Let $H$ be as in Theorem A. Then for any $k>0$, there exists a unique positive solution of \eqref{A1} with $\mu=k\gd_0$  (where $\gd_0$ is the Dirac mass  at the origin). This solution is denoted by $u^\Gw_{k,0}$.

Furthermore,
\bel{weaksing} u^{\Gw}_{k,0}(x)=kP^\Gw(x,0)(1+o(1)) \q \text{as } x \to 0. \ee
and there exists $d_k>0$ such that
\bel{com} d_k P^\Gw(x,0)<u^\Gw_{k,0}(x) < k P^\Gw(x,0) \forevery x \in \Gw. \ee
}

 Obviously, $u^\Gw_{k,0}$ is a moderate solution and is called a \emph{weakly singular solution}.  It follows from \eqref{weaksing} that the sequence $\{u^\Gw_{k,0}\}$ is increasing. Moreover, this sequence is uniformly bounded in any compact subset of $\Gw$. Therefore
 $$u^\Gw_{\infty,0}:=\lim u^\Gw_{k,0}$$
 is a solution of \req{A0}. Clearly this solution is not moderate.

 When there is no danger of confusion we drop the upper index writing simply $u_{k,0}$ and $u_{\infty,0}$.

 Denote by $\CU_0^\Gw$ the family of positive \emph{non-moderate} solutions of \req{A0} \sth
 $u\in C(\ovl \Gw\sms \{0\})$ and $u=0\;\text{on}\; \bdw\sms \{0\}$. If $u$ is such a solution we say that it is \emph{a strongly singular solution}. In  the next two theorems we consider solutions of this type.
    \medskip

\noindent{\bf Theorem D.} {\it Under the assumptions of theorem C, $u^\Gw_{\infty,0}\in \CU_0^\Gw$.  Furthermore $u^\Gw_{\infty,0}$
 is the minimal element of $\CU^\Gw_0$.} \medskip

Let $S^{N-1}$ be the unit sphere, $\BBR^N_+=[x_N >0]$, $S_+^{N-1}=S^{N-1} \cap \BBR^N_+$ the upper hemisphere and $(r,\gs) \in \BBR_+ \ti S^{N-1}$ the spherical coordinates in $\BBR^N $. Denote by $\nabla' $ and  $\Gd'$ the covariant derivative on $S^{N-1}$ identified with the tangential derivative and the Laplace-Beltrami operator on $S^{N-1}$ respectively.

 As mentioned before, if $H$ is as in \eqref{multi},  \eqref{A0} admits a similarity transformation. However, there is no similarity transformation when $H$ satisfies \eqref{add} unless $p = \frac{q}{2-q}$. When  $p \neq \frac{q}{2-q}$ there is competition between $u^p$ and $|\nabla u|^q$.
 When $p>\frac{q}{2-q}$ the dominant term is $u^p$;  when $p<\frac{q}{2-q}$ the dominant term is $|\nabla u|^q$.
This fact is reflected in the next theorem.

We assume that the set of coordinates is placed so that $0\in \bdw$, $x_N=0$ is tangent to $\bdw$ at $0$ and the positive $x_N$ axis points into the domain.
\medskip

\noindent{\bf Theorem E.}\hskip 2mm
\textit{Assume that either $H$ satisfies  \req{multi}, $0<N(p+q-1)<p+1$ and $p\geq1$ or $H$ satisfies \req{add} and $m_{p,q}<p_c$ where $p_c$ and $m_{p,q}$ are given by \eqref{pc} and \eqref{mpq} respectively. Then:}

(i)\textit{ $\CU^{\Gw}_0$ consists of a single element $u^{\Gw}_{\infty,0}$. In other words, $u^{\Gw}_{\infty,0}$ is the unique strongly singular solution of \req{A0} with singularity at $0$.}
\medskip

(ii)\textit{  Put $r=|x|$, $\gs=\frac{x}{r}$. Then
\begin{equation}\label{U0infty}
 \lim_{x \in \Gw, r \to 0}r^\gb u^{\Gw}_{\infty,0}(x) = \gw(\gs)
\end{equation}
locally uniformly on $S^{N-1}_+$ where
\bel{beta1}\gb=\gb_1:= \frac{2-q}{p+q-1} \q\text{if $H$ satisfies \req{multi}},\ee
and
\bel{beta2} \gb=\gb_2:= \frac{2}{m_{p,q}-1} \q\text{if $H$ satisfies \req{add}}.\ee}

\textit{The function $\gw$ is the unique solution of the problem
\begin{equation}\label{Hgw} \left\{\BAL
 -\Gd'\gw+F(\gw,\nabla'\gw)&=0 \q\text{in }S^{N-1}\\
 \gw&=0 \q\text{on }\prt S^{N-1}
\EAL \right. \end{equation}
where $F=F_i(s,\gx)$  $(i=1,\ldots,4)$,   $(s,\gx)\in \BBR_+\ti S^{N-1}$, is given by,
\begin{equation}\label{i=1...4}\BAL
 &F_1(s,\gx)=s^p(\gb_1^2\, s^2+\abs{\gx}^2)^{\frac{q}{2}}-\gb_1(\gb_1+2-N)s,&& \text{when $H$ satisfies }\req{multi}\\
 &\text{and, when $H$ satisfies \req{add},}&&\\
 &F_2(s,\gx)= s^p + (\gb_{2}^2\, s^2+\abs{\gx}^2)^{\frac{q}{2}}-\gb_{2}(\gb_{2}+2-N)s, && \text{if } p=\frac{q}{2-q}\\
 &F_3(s,\gx)= s^p-\gb_{2}(\gb_{2}+2-N)s,  && \text{if } p>\frac{q}{2-q}\\
 &F_4(s,\gx)=(\gb_{2}^2\, s^2+\abs{\gx}^2)^{\frac{q}{2}}-\gb_{2}(\gb_{2}+2-N)s.  && \text{if } p<\frac{q}{2-q}\\.
\EAL\end{equation}
}
\textit{The unique solution of \req{Hgw} with $F=F_i$ will be denoted by $\gw_i$, $i=1,\ldots,4$. When $i=1,2$ we actually have $u^{\BBR^N_+}_{\infty,0}(x)= r^{-\gb_i}\gw_i(\gs)$.}
\medskip

\Remark We note that $r^{-\gb_2}\gw_3$ (resp. $r^{-\gb_2}\gw_4$) behaves near the origin like the corresponding strongly singular solution of $-\Gd u+u^p=0$ (resp. $-\Gd u+|\nabla u|^q=0$).
\vskip 2mm

Next we present a removability result which implies that the conditions on $p,q$ for $H$ to be subcritical are sharp.
\medskip

\noindent{\bf Theorem F} {\it Assume that $H$ satisfies either \eqref{multi} with $N(p+q-1) \geq p+1$
or \eqref{add} with $m_{p,q} \geq p_c$. If $u\in C(\overline\Gw\setminus\{0\})\cap C^2(\Gw)$ is a nonnegative solution of \eqref{A0} vanishing on $\prt\Gw\setminus\{0\}$ then $u\equiv 0$.} \medskip

When $H$ satisfies \eqref{add}, the removabilty result  is based on the corresponding results for  \req{A0} with $H=u^p$ and  $H=|\nabla u|^q$.

When $H$ satisfies \eqref{multi} and $N(p+q-1)>p+1$ we use a similarity transformation to show that there is no solution with isolated singularity. The case $N(p+q-1)=p+1$ is a bit more delicate. We first establish the removability result for the half-space and use it, together with some regularity results up to the boundary (see \cite{Lib}) to derive the result in bounded domains of class $C^2$.

If $q=2$, one can obtain removability result by a change of unknown.
\medskip

\Remark Theorems C, F and G provide a \emph{complete characterization of the positive solutions $u$ of \req{A0} in $\Gw$ \sth $u=0$ on $\bdw$ except at one point.}
\medskip

In the case where $H$ satisfies \req{add} we also consider solutions that blow-up strongly on an arbitrary compact set $K\sbs \bdw$.
\medskip

\noindent{\bf Theorem G} \textit{Assume $H$ satisfies \eqref{add} and $m_{p,q} \leq p_c$ where $m_{p,q}$ and $p_c$ are given by \eqref{mpq} and \eqref{pc} respectively. Let $K$ be a compact subset of $\prt \Gw$. Denote by $\CU_K$ the family of all positive solutions $u$ of \eqref{A0} such that $\CS(u)=K$ (see \rdef{bdtrace}) and $u=0$ on $\Gw \sms K$. Then there exist a minimal element $u_K$ and a maximal element $U_K$ of $\CU_K$ in the sense that $u_K \leq u \leq U_K$ for every $u \in \CU_K$. Morover, for every $y \in K$ and $\gg \in (0,1)$, there exist $r$ depending on $\gg$ and $C$ depending on $N$, $p$, $q$, $\gg$ and the $C^2$ characteristic of $\Gw$ such that}
\bel{2K} U_K(x) \leq C u_K(x) \forevery x \in C_{\gg,r}(y):=\{x \in \Gw: \rho(x) \geq \gg |x-y|\} \cap B_r(y). \ee

This extends a result of \cite{BaGi} on existence of large solutions.

\medskip

The paper is organized as follows. In section 2, we establish some estimates on positive solution of \eqref{A0} and its gradient, and recall some estimates concerning weak $L^p$ space.  Section 3 is devoted to the proof of Theorems A, B and various results on {\it boundary trace}. In section 4, we provide a complete description of isolated singularities on the boundary (Theorems C,D,E). Boundary value problem with unbounded measure data for \eqref{A0} and $H$ as in \eqref{add} is discussed in Section 5 (Theorem G). In section 6, we demonstrate the removability result in the supercritical case (Theorem F). In the appendix, a uniqueness result for a class of quasilinear elliptic equations is proved.

\section{Preliminaries}
Throughout the present paper, we denote by $c$, $c_1$, $c_2$, $C$,...positive constants which may vary from line to line. If necessary the dependence of these constants will be made precise.
The following comparison principle can be found in \cite[Theorem 9.2]{GT}.
\bprop{comparison} Assume $H: D \ti \BBR_+ \ti \BBR^N \to \BBR_+$ is nondecreasing with respect to $u$ for any $(x,\gx) \in D \ti \BBR^N$, continuously differentiable with respect to $\gx$ and $H(x,0,0)=0$. Let $u_1$, $ u_2 \in C^2(D) \cap C(\ovl D)$ be two nonnegative solutions of \eqref{A0}. If
$$ -\Gd u_1 + H \circ u_1 \leq -\Gd u_2 + H \circ u_2 \quad \text{in } D $$
and $u_1 \leq u_2$ on $\prt D$. Then $u_1 \leq u_2$ in $D$.
\es

Next, for $\gd>0$, we set
	$$ \Gw_\gd=\{x \in \Gw: \rho(x) < \gd\},\q D_\gd=\{x \in \Gw: \rho(x) > \gd\}, \q \Gs_\gd=\{x \in \Gw: \rho(x) = \gd\}. $$
	
\bprop{delta0} There exists $\gd_0>0$ such that

(i) For every point $x \in \ovl \Gw_{\gd_0}$, there exists a unique
point $\gs_x \in \prt \Gw$ such that $|x -\gs_x|=\rho(x)$. This implies
$x=\gs_x - \rho(x){\bf n}_{\gs_x}$.

(ii) The mappings $x \mapsto \rho(x)$ and $x \mapsto \gs_x$ belong to
$C^2(\ovl \Gw_{\gd_0})$ and $C^1(\ovl \Gw_{\gd_0})$ respectively. Furthermore,
$\lim_{x \to \gs_x}\nabla \rho(x) = - {\bf n}_{\gs_x}$.

\es
In the sequel, we can assume that $\gd_0< \norm{\Gd \rho}^{-1}_{L^\infty(\Gw)}$. The next results provide a-priori estimates on positive solutions and their gradient.
\bprop{est-multi} Assume $H$ satisfies \eqref{multi} with $p \geq 0$, $0\leq q <2$, $p+q>1$. Let $u \in C^2(\Gw)$ be a positive solution of equation \eqref{A0}. Then
\bel{est2} u(x) \leq  \Gl_1\rho(x)^{-\gb_1} + \Gl_1'\norm{u}_{L^1(D_{\frac{\gd_0}{2}})} \quad \forall x \in \Gw,\ee
\bel{grad} \abs{\nabla u(x)} \leq \tl \Gl_1\,\rho(x)^{-\gb_1-1} \forevery x \in \Gw\ee
where $\gb_1$ is defined in \eqref{beta1}, $\tl \Gl_1=\tl \Gl_1(N,p,q,\gd_0,\norm{u}_{L^1(D_{\frac{\gd_0}{2}})})$,  $\Gl_1'=\Gl_1'(N,\gd_0)$,  and
\bel{Lambda1}  \Gl_1=\left(\myfrac{\gb_1+2}{ \gb_1^{q-1}}\right)^{\frac{1}{p+q-1}}. \ee
\es
\Proof Put $ M_{\gd_0}=\max\{u(x): x \in \ovl D_{\gd_0}\}$. For each $\gd \in (0,\gd_0)$, we set
$$ w_{\gd}(x)=\Gl_1(\rho(x)-\gd)^{-\gb_1}+M_{\gd_0} \q x \in D_\gd.$$
By a simple computation, we obtain
$$ -\Gd w_{\gd} +  w_{\gd}^p\abs{\nabla w_{\gd}}^q > 0 \q \text{in } \Gw_{\gd_0}\sms \ovl \Gw_{\gd}.$$
Since $w_\gd \geq u$ on $\Gs_{\gd} \cup \Gs_{\gd_0}$, by the comparison principle \rprop{comparison}, $u \leq w_\gd$ in $\Gw_{\gd_0} \sms \ovl \Gw_{\gd}$. Letting $\gd \to 0$ yields
\bel{est1} u(x) \leq \Gl_1 \rho(x)^{-\gb_1} + M_{\gd_0} \quad \forall x \in \Gw. \ee
Since $u$ is subharmonic, by \cite[Theorem 1]{Tr1}, there exists $\Gl_1'=\Gl_1'(N,\gd_0)$ such that $\Gl_1' \norm{u}_{L^1(D_{\gd_0/2})}>M_{\gd_0}$. This, along with \eqref{est1}, implies \eqref{est2}.   \medskip

Next we prove \eqref{grad}. Fix $x_0 \in \Gw$ and set $d_0=\frac{1}{3}\rho(x_0)$, $y_0=\frac{1}{d_0}x_0$ and
$$ \quad M_{0}=\max\{u(x): x \in B_{2d_0}(x_0)\}, \q
u_0(y)=\myfrac{u(x)}{M_0},\quad  y=\frac{1}{d_0}x \in B_2(y_0). $$
Then $\max\{u_0(y): y \in B_2(y_0)\}=1$ and $ -\Gd u_0 +   M_0^{p+q-1}d_0^{2-q}u_0^p|\nabla u_0|^q = 0$ in $B_2(y_0)$. By \cite{La},  there exists a positive constant $c=c(N,p,q,\gd_0,\norm{u}_{L^1(D_{\frac{\gd_0}{2}})})$ such that $\max_{B_1(y_0)}|\nabla u_0| \leq c$. Consequently,
	$$ \max_{B_{d_0}(x_0)}|\nabla u| \leq cd_0^{-1}\max_{B_{2d_0}(x_0)}u $$
which implies \eqref{grad}. \qed

\bprop{est-add} Assume $H$ satisfies \eqref{add} with $p > 1$, $1< q <2$. Let $u \in C^2(\Gw)$ be a positive solution of equation \eqref{A0}. Then
\bel{est2-add} u(x) \leq  \Gl_2\rho(x)^{-\gb_2}  \ee
\bel{grad-add} \abs{\nabla u(x)} \leq \tl \Gl_2\,\rho(x)^{-\gb_2-1} \forevery x \in \Gw\ee
where $\gb_2$ is defined in \eqref{beta2}, $\Gl_2=\Gl_2(N,p, q,\gd_0)$ and $\tl \Gl_2=\tl \Gl_2(N,p,q,\gd_0)$.
\es
\Proof Since $u$ is a subsolution of \eqref{Ap}, it follows from Keller-Osserman \cite{MVbook} that
\bel{e1-add} u(x) \leq c\rho(x)^{-\frac{2}{p-1}} \forevery x \in \Gw \ee
where $c=c(N,p)$. By a similar argument as in the proof of \rprop{est-multi}, we deduce that
\bel{e2-add} u(x) \leq  c'\rho(x)^{-\gb_2} + c''\norm{u}_{L^1(D_{\frac{\gd_0}{2}})} \forevery x \in \Gw \ee
where $c'=c'(p,q)$ and $c''=c''(N,\gd_0)$. Combining \eqref{e1-add} and \eqref{e2-add} implies \eqref{est2-add}. Finally, we derive \eqref{grad-add} from \rprop{est-multi} as in the proof of \rprop{est-multi}. \qed \medskip

Denote  by $G^\Gw$ (resp. $P^\Gw$)  the Green kernel (resp. the Poisson kernel) in $\Gw$, with corresponding operators $\BBG^\Gw$ (resp. $\BBP^\Gw$). We denote by $\GTM_{\rho^\ga}(\Gw)$, $\ga \in [0,1]$, the space of Radon measures $\gm$ on $\Gw$ satisfying $\int_\Gw\rho^\ga d|\gm|<\infty$, by $\GTM(\prt \Gw)$ the space of bounded Radon measures on $\prt \Gw$ and by $\GTM^+(\prt \Gw)$ the positive cone of $\GTM(\prt \Gw)$.

Denote  $L^p_w(\Gw;\tau)$, $1 \leq p < \infty$, $\tau \in \GTM^+(\Gw)$, the weak $L^p$ space defined as follows: a measurable function $f$ in $\Gw$ belongs to this space if there exists a constant $c$ such that
\bel{distri} \gl_f(a;\tau):=\tau(\{x \in \Gw: |f(x)|>a\}) \leq ca^{-p}, \forevery a>0. \ee
The function $\gl_f$ is called the distribution function of $f$ (relative to $\tau$). For $p \geq 1$, denote
$$ L^p_w(\Gw;\tau)=\{ f \text{ Borel measurable}: \sup_{a>0}a^p\gl_f(a;\tau)<\infty\} $$
and
\bel{semi} \norm{f}^*_{L^p_w(\Gw;\tau)}=(\sup_{a>0}a^p\gl_f(a;\tau))^{\frac{1}{p}}. \ee
The $\norm{.}_{L^p_w(\Gw;\tau)}$ is not a norm, but for $p>1$, it is equivalent to the norm
\bel{normLw} \norm{f}_{L^p_w(\Gw;\tau)}=\sup\left\{ \frac{\int_{\gw}|f|d\tau}{\tau(\gw)^{1/p'}}:\gw \sbs \Gw, \gw \text{ measurable }, 0<\tau(\gw)<\infty \right\}. \ee
More precisely,
\bel{equinorm} \norm{f}^*_{L^p_w(\Gw;\tau)} \leq \norm{f}_{L^p_w(\Gw;\tau)} \leq \myfrac{p}{p-1}\norm{f}^*_{L^p_w(\Gw;\tau)} \ee
Notice that, for every $\ga>-1$,
$$L_w^p(\Gw;\rho^{\ga}dx) \sbs L_{\rho^{\ga}}^{s}(\Gw) \forevery s \in [1,p). $$

The following useful estimates involving Green and Poisson operators can be found in \cite{BVi} (see also \cite{MVbook} and \cite{V2}).
\bprop{P1} For any $\ga\in [0,1]$, there exist a positive constant $c_1$ depending on $\ga$,  $\Gw$ and $N$ such that
\bel{E1-1} \BA{lll}
\norm{\BBG^{\Gw}[\gn]}_{L_w^{\frac{N+\ga}{N+\ga-2}}(\Gw;\rho^\ga dx)}+\norm{\nabla\BBG^{\Gw}[\gn]}_{L_w^{\frac{N+\ga}{N+\ga-1}}(\Gw;\rho^\ga dx)}
\leq c_1\norm \gn_{\mathfrak M_{\rho^\ga}(\Gw)},
\EA \ee
\bel{E1-3} \BA{lll}
\norm{\BBP^{\Gw}[\gm]}_{L_w^{\frac{N+\ga}{N-1}}(\Gw,\rho^\ga dx)}+\norm{\nabla\BBP^{\Gw}[\gm]}_{L_w^{\frac{N+1}{N}}(\Gw;\rho^\ga dx)}
\leq c_1\norm \gm_{\GTM(\prt \Gw)},
\EA \ee
for any $\gn\in \GTM_{\rho^\ga}(\Gw)$ and any $\gm\in \GTM(\prt\Gw)$ where
$$ \norm \gn_{\mathfrak M_{\rho\ga}(\Gw)}:=\myint{\Gw}{}\rho^\ga d|\gn| \q \text{and} \q \norm \gm_{\GTM(\prt \Gw)}=\myint{\prt \Gw}{}d|\gm|. $$
\es\medskip

\section{Boundary value problem with measures and boundary trace}
\subsection{The Dirichlet problem}
\noindent{\bf Proof of Theorem A.} We deal with the case when $H$ satisfies \eqref{multi}. The case $H$ satisfies \eqref{add} is simpler and can be treated in a similar way.

Let $\{\gm_n\}$ be a sequence of  positive functions in $C^1(\prt \Gw)$  converging weakly to $\gm$. There exists a positive constant $c_2$ independent of $n$ such that $\norm{\gm_n}_{L^1(\prt \Gw)} \leq c_2\norm{\gm}_{\GTM(\prt \Gw)}$ for all $n$. Consider the following problem
	\bel{v_n} \left\{ \BA{lll}
		- \Gd v +  (v+ \BBP^\Gw[\gm_n])^p|\nabla (v + \BBP^{\Gw}[\gm_n])|^q &= 0 \qq &\text{ in } \Gw\\
		\phantom{- \Gd  +g(\abs{\nabla (v + \BBP^{\Gw}[\gm_n])}),,,,---}
		v &= 0  &\text{ on } \prt \Gw.
	\EA \right. \ee
Since $0$ and $-\BBP^{\Gw}[\gm_n]$ are respectively supersolution and subsolution of \eqref{v_n}, by \cite[Theorem 6.5]{KaKr} there exists a solution $v_n \in W^{2,s}(\Gw)$ with $1<s<\ity$ to problem \eqref{v_n} satisfying $-\BBP^{\Gw}[\gm_n] \leq v_n \leq 0$. Thus $u_n = v_n + \BBP^{\Gw}[\gm_n]$ is a solution of
	\bel{u_n} \left\{ \BA{lll}
		- \Gd u_n +   u_n^p\abs{\nabla u_n}^q &= 0 \qq &\text{ in } \Gw\\
		\phantom{- \Gd  + g(\abs{\nabla u_n}),}
		u_n &= \gm_n  &\text{ on } \prt \Gw.
	\EA \right. \ee
By the maximum principle, such solution is the unique solution of \eqref{u_n}. \medskip

\noindent{\it Assertion 1:} $\{u_n\}$ and $\{\abs{\nabla u_n}\}$ remain uniformly bounded respectively in $ L_w^\frac{N}{N-1}(\Gw)$ and $L_w^\frac{N+1}{N}(\Gw;\rho dx)$.

Let $\gx$ be the solution to
	 \bel{eta} - \Gd \gx = 1 \text{ in } \Gw, \q \gx = 0  \text{ on } \prt \Gw, \ee
then there exists a constant $c_3>0$ such that $c_3^{-1}<-\frac{\prt \gx}{\prt \bf n} < c_3$ on $\prt \Gw$ and $c_3^{-1}\rho \leq \gx \leq c_3\rho$ in $\Gw$.
By multiplying the equation in \eqref{u_n} by $\gx$ and integrating on $\Gw$, we obtain
	\bel{EM2-1} \myint{\Gw}{}u_n dx +  \myint{\Gw}{}u_n^p\abs{\nabla u_n}^q\rho dx \leq c_4\norm{\gm}_{\GTM(\prt \Gw)} \ee
where $c_4$ is a positive constant independent of $n$. From \rprop{P1} and by noticing that $u_n \leq \BBP^{\Gw}[\gm_n]$, for every $\ga \in [0,1]$, we get
	\bel{EM2-2} \norm{u_n}_{L_w^\frac{N+\ga}{N-1}(\Gw;\rho^\ga dx)} \leq \norm{\BBP^\Gw[\gm_n]}_{L_w^{\frac{N+\ga}{N-1}}(\Gw;\rho^\ga dx)} \leq c_1\norm{\gm_n}_{L^1(\prt \Gw)} \leq c_1c_2\norm{\gm}_{\GTM(\prt \Gw)}. \ee
Again, from \rprop{P1} and  \eqref{EM2-1}, we derive that
	\bel{EM2-3} \norm{\nabla u_n}_{L_w^\frac{N+1}{N}(\Gw;\rho dx)} \leq c_1\left(\norm{u_n^p\abs{\nabla u_n}^q}_{L_\rho^1(\Gw)} + \norm{\gm_n}_{L^1(\prt \Gw)}\right)  \leq c_5\norm{\gm}_{\GTM(\prt \Gw)} \ee
where $c_5$ is a positive constant depending only on $\Gw$ and $N$.
Thus Assertion 1 follows from \eqref{EM2-2} and \eqref{EM2-3}.
	
By regularity results for elliptic equations \cite{Mi}, there exist a subsequence, still denoted by $\{u_n\}$, and a function $u$ such that $\{u_n\}$ and $\{|\nabla u_n|\}$ converges to $u$ and $|\nabla u|$  a.e. in $\Gw$.\medskip

\noindent {\it Assertion 2:} $\{u_n\}$ converges to $u$ in $L^1(\Gw)$.

Indeed, by taking $\ga=0$ in \eqref{EM2-2}, we derive $\{u_n\}$ is uniformly bounded in $L_w^{\frac{N}{N-1}}(\Gw)$. Therefore, $\{u_n\}$ is uniformly bounded in $L^r(\Gw)$ for any $r \in [1,\frac{N}{N-1})$. By Holder inequality, $\{ u_n\}$ is uniformly integrable in $L^1(\Gw)$.  Thus Assertion 2 follows from Vitali's convergence theorem.\medskip

\noindent {\it Assertion 3:} $\{u_n^p|\nabla u_n|^q\}$ converges to $u^p|\nabla u|^q$ in $L_\rho^1(\Gw)$.

Indeed, by taking $\ga=1$ in \eqref{EM2-2}, one derives that $\{u_n\}$ is uniformly bounded in $L_w^{\frac{N+1}{N-1}}(\Gw;\rho dx)$. Therefore, $\{u_n\}$ is uniformly bounded in $L_\rho^r(\Gw)$ for every $r \in [1,\frac{N+1}{N-1})$. By \eqref{EM2-3}, $\{|\nabla u_n|\}$ is uniformly bounded in $L_\rho^s(\Gw)$ for every $s \in [1,\frac{N+1}{N})$. Since $N(p+q-1)<p+1$, we can choose $r$ and $s$ close to $\frac{N+1}{N-1}$ and $\frac{N+1}{N}$ respectively so that $\frac{p}{r}+\frac{q}{s}<1$. By Holder inequality, $\{u_n^p|\nabla u_n|^q\}$ is uniformly integrable in $L_\rho^1(\Gw)$. Thus Assertion 3 follows from Vitali's convergence theorem. \medskip

For every $\zeta \in C_0^2(\ovl \Gw)$, we have
	\bel{EM2-8} \myint{\Gw}{}(-u_n\Gd \zeta +  u_n^p \abs{\nabla u_n}^q\zeta)dx=-\myint{\prt \Gw}{}\gm_n\myfrac{\prt \zeta}{\prt \bf n}dS.
	\ee
Due to Assertions 2 and 3, by letting $n \to \infty$ in \eqref{EM2-8} we obtain \eqref{A2}; so $u$ is a solution of \eqref{A1}. By \rprop{P1}, $u \in L_w^\frac{N}{N-1}(\Gw)$ and  $\abs{\nabla u} \in L_w^\frac{N+1}{N}(\Gw;\rho dx)$.

Next, let $\{\gm_n\}$ be a sequence of positive finite measures on $\prt\Gw$ which converges weakly to a positive finite measure $\gm$ and $\{u_{\gm_n}\}$ be a sequence of corresponding solutions of \eqref {u_n}. Then by using a similar argument as in Assertions 2 and 3, we deduce that there exists a subsequence such that $\{ u_{\gm_{n_k}}\}$ converges to a solution $ u_\gm$ of \eqref {A1} in $L^1(\Gw)$ and $\{H \circ u_{\gm_{n_k}}\}$ converges to $H \circ u$ in $L^1_\rho(\Gw)$. \qed \medskip

Using Theorem A one can establish a slightly stronger type of stability.
\bcor{var-stab} Let $H$ be as in theorem A.  Let $\{a_n\}$ be a decreasing sequence converging to $0$, $\gm$ be a bounded positive measure on $\prt \Gw$ and $\{\gm_n\}$ be a sequence of bounded positive measure on $\Gs_{a_n}$ converging weakly to $\gm$. Let $\{u_{\gm_n}\}$ be a sequence of corresponding solutions of \eqref {u_n} in $D_{a_n}$. Then there exists a subsequence such that $\{ u_{\gm_{n_k}}\}$ converges in $L^1(\Gw)$ to a solution $ u_\gm$ of \eqref {A1} and $\{H \circ u_{\gm_{n_k}}\}$ converges to $H \circ u$ in $L^1_\rho(\Gw)$.
\es
\Proof As above, we consider the case $H$ satisfies \eqref{multi} because the case $H$ satisfies \eqref{add} can be proved by a similar argument. We extend $ u_{\gm_n}$ and $|\nabla u_{\gm_n}|$ by zero outside $\ovl D_{a_n}$ and still denote them by the same expressions. By regularity results for elliptic equations \cite{Mi}, there exist a subsequence, still denoted by $\{u_{\gm_n}\}$, and a function $u$ such that $\{u_{\gm_n}\}$ and $\{|\nabla u_{\gm_n}|\}$ converges to $u$ and $|\nabla u|$  a.e. in $\Gw$. Let $G\subset\Gw$ be a Borel set and put $G_n=G\cap D_{a_n}$. By using similar argument as in Assertion 2 in the proof of theorem A, due to the estimate $||\BBP^{\Gw}[\gm]|_{_{\Gs_{a_n}}}||_{L^1(\Gs_{a_n})}\leq c_7\norm{\gm}_{\mathfrak M(\Gs)}$, we derive
\bel{EM2-5+1}\BA {lll} \myint{G_n}{}u_{\gm_n} dx&\leq |G_n|^{\frac{1}{N}}\norm{u_{\gm_n}}_{L_w^{\frac{N}{N-1}}(D_{a_n})}
\leq c_1c_2|G_n|^{\frac{1}{N}}\norm{\BBP^\Gw[\gm]|_{_{\Gs_{a_n}}}}_{L^1(\Gs_{a_n})}\\[4mm]\phantom{\myint{G_n}{}u_\gd dx}
&\leq c_1c_2c_7|G|^{\frac{1}{N}}\norm{\gm}_{\GTM(\Gs)}.
\EA\ee
Hence $\{u_{\gm_n}\}$ is uniformly integrable. Therefore due to Vitali's convergence theorem, up to a subsequence, $\{u_{\gm_n}\}$ converges to $u$ in $L^1(\Gw)$.

Set $\rho_n(x):=(\rho(x)-a_n)_+$. By proceeding as in Assertion 3 of the proof of Theorem A and notice that $\int_{G_n}^{}\rho_n dx \leq \int_{G}^{}\rho dx$, we derive that $\{u_{\gm_n}^p|\nabla u_{\gm_n}|^q\}$ is uniformly integrable. Therefore by Vitali's convergence, up to a subsequence, $\{ u_{\gm_n}^p |\nabla u_{\gm_n}|^q\}$ converges to $u^p|\nabla u|^q$ in $L^1_\rho(\Gw)$.

Finally, if $\gz\in C_0^2(\ovl \Gw)$ we denote by $\gz_n$ the solution of
\bel{EM2-9} -\Gd \gz_n=-\Gd\gz \text{ in } D_{a_n}, \q \gz_n=0 \text{ on } \prt D_{a_n}. \ee
Then $\gz_n \in C_0^2(\ovl \Gw_{a_n})$, $\gz_n \to \gz$ in $C^2(\Gw)$ and $\sup_n\norm{\gz_n}_{C^2(\ovl \Gw_{a_n})}<\infty$. Since
	\bel{EM2-10} \myint{D_{a_n}}{}(- u_{\gm_n}\Gd \zeta_n +   u_{\gm_n}^p\abs{\nabla u_{\gm_n}}^q\zeta_n)dx=-\myint{\Gs_{a_n}}{}\myfrac{\prt \zeta_n}{\prt \bf n}d \gm_n,
	\ee
by letting $n \to \infty$, we deduce that $u$ is a solution of \eqref{A1}. \qed \medskip

\Remark  Let $\gm \in \GTM^+(\prt \Gw)$. It follows from \rprop{est-multi} and \rprop{est-add} that there exists a constant $c$ depending on $N$, $p$, $q$, $\Gw$ and $\norm{\gm}_{\GTM(\prt \Gw)}$ such that for every positive solution $u$ of \eqref{A1} there holds
	\bel{est1*} u(x) \leq c\rho(x)^{-\gb_i} \forevery x \in \Gw, \ee
	\bel{est2*}|\nabla u(x)| \leq c\rho(x)^{-\gb_i-1} \forevery x \in \Gw \ee
	where
	$$  i=\left\{ \BA{lll} 1 \text{ if } H \text{ satisfies } (\ref{multi}) \\
	2 \text{ if } H \text{ satisfies } (\ref{add}).
	\EA \right. $$
	
Using these facts and \rcor{var-stab} we obtain the following monotonicity result.
\bcor{max} Let $H$ be as in theorem A. For any $\gm \in \GTM^+(\prt \Gw)$, there exists a maximal solution $U_\gm$ of \eqref{A1}. Moreover, if $\gm, \gn \in \GTM^+(\prt \Gw)$ such that $\gm \leq \gn$ then $U_\gm \leq U_\gn$. \es
\Proof For each $\gd>0$, let $U:=U_{\gm,\gd}$ be the solution of
\bel{U} \left\{ \BA{lll} -\Gd U + H \circ U &= 0 \quad &\text{in } D_\gd \\
\phantom{-\Gd U + H \circ,}
U &= \BBP^\Gw[\gm] &\text{on }  \Gs_\gd.
\EA \right. \ee
By the comparison principle, $0 \leq U_{\gm,\gd} \leq \BBP^\Gw[\gm]$, hence $\{ U_{\gm,\gd}\}$ is decreasing as $\gd \to 0$. Put $U_\gm:=\lim_{\gd \to 0}U_{\gm,\gd}$ then by \rcor{var-stab} $U_\gm$ is a solution of \eqref{A1}. If $u$ is a positive solution of \eqref{A1} then by the comparison principle $0 \leq u \leq \BBP^\Gw[\gm]$ in $\Gw$. Therefore $u \leq U_{\gm,\gd}$ in $D_\gd$ for every $\gd>0$. Letting $\gd \to 0$ implies $u \leq U_\gm$.

Next, if $\gm \leq \gn$ then $\BBP^\Gw[\gm] \leq \BBP^\Gw[\gn]$. Hence $U_{\gm,\gd} \leq U_{\gn,\gd}$ for every $\gd>0$ and therefore $U_\gm \leq U_\gn$. \qed \medskip

\noindent{\bf Proof of Theorem B.} The strategy is the same as in the proof of \rth{uniqueA} so we only sketch the main technical modifications. Let $u$ be a positive solution of \eqref{A1} then $u \leq U_\gm$. Let $\{\gm_n\}$ be a sequence of functions in $C^1(\prt \Gw)$ converging weakly to $\gm$. For $k>0$, denote by $T_k$ the truncation function, i.e. $T_k(s)=\max(-k,\min(s,k))$. For every $n>0$, denote by $u_n$ and $U_{\mu,n}$ respectively the solutions of
\bel{o2} - \Gd u_{n} + T_n(H \circ u) = 0 \quad \text{in } \Gw, \quad u_{n}=\mu_n \quad \text{on } \prt \Gw.
\ee
\bel{o2'} - \Gd U_{\mu,n} + T_n(H \circ U_\mu) = 0 \quad \text{in } \Gw, \quad U_{\mu,n}=\mu_n \quad \text{on } \prt \Gw.
\ee
By local regularity theory for elliptic equations (see, e.g., \cite{Mi}), $u_{n} \to u$ and $U_{\mu,n} \to U_\mu$ in $C^1_{loc}(\Gw)$.
From \eqref{o2} and \eqref{o2'} we obtain
\bel{o6} \left\{ \BA{lll} - \Gd (U_{\mu,n}-u_{n}) &= -T_n(H \circ U_\mu) + T_n(H \circ u) \quad \text{in } \Gw \\
\phantom{--,,}
U_{\mu,n}-u_{n} &= 0 \quad \text{on } \prt \Gw. \EA \right. \ee
We shall prove that $U_\mu = u$. By contradiction, we assume that $M:=\sup_\Gw(U_\mu-u) \in (0,\infty]$. Let $0<k<M$. From \eqref{o6}, Kato's inequality \cite{MVbook} and the fact that $u \leq U_\mu$, we get
\bel{o7} \BA{lll} -\Delta(U_{\mu,n}-u_{n}-k)_+ &\leq (T_n(H \circ u) - T_n(H \circ U_{\mu,n}))\chi_{_{E_{n,k}}} \\
&\leq ||\nabla U_\mu|^q-|\nabla u|^q|\chi_{_{E_{n,k}}}
\EA \ee
where $E_{n,k}=\{x \in \Gw: u_{1,n}-u_{2,n}>k\}$. We next proceed as in the proof of \rth{uniqueA} in order to get a contradiction. Thus $u= U_\gm$. \qed \medskip

As a consequence, we obtain the following comparison principle
\bcor{comparison} Under the assumption of Theorem B, if $u_1$ and $u_2$ be respectively positive sub and supersolution solution of \eqref{A0} such that $\tr(u_1) \leq \tr(u_2)$ then $u_1 \leq u_2$ in $\Gw$.
\es
\Proof We first observe that if $u_1$ and $u_2$ are both solution of \eqref{A0} then by Theorem A and B, $u_1 \leq u_2$.

Next we consider the case $u_1$ and $u_2$ are respectively sub and super solution. For $\gd>0$, let $v_{i,\gd}$, $i=1,2$ be the solution of
\bel{vi}  \left\{ \BAL -\Gd v + H \circ v &= 0 \q \text{in } D_\gd \\
v &= u_i \q \text{on } \Gs_\gd
\EAL \right.  \ee
By the comparison principle, $u_1 \leq v_{1,\gd}$ and $v_{2,\gd} \leq u_2$ in $D_\gd$. Therefore $\{v_{1,\gd}\}$ and $\{v_{2,\gd}\}$ are respectively increasing and decreasing as $\gd \to 0$. By \rcor{var-stab} and Theorem B, $v_i:=\lim_{\gd \to 0}v_{i,\gd}$ is the solution of \eqref{A0} with boundary trace $\mu_i$. Moreover, $u_1 \leq v_1$ and $v_2 \leq u_2$. Since $\mu_1 \leq \mu_2$, by the above observation, $v_1 \leq v_2$. Thus $u_1 \leq v_1 \leq v_2 \leq u_2$. \qed \medskip

When $H$ satisfies \eqref{multi}, the question of uniqueness remains open, but we can show that any positive solution of \eqref{A1} behaves like $ U_\gm$ near the boundary. Before stating the result, we need the following definition
\bdef{potential} A nonnegative superharmonic function is called a \emph{potential} if its largest harmonic minorant is zero.
\es
\bprop{nontan} Let $\gm \in \GTM^+(\prt \Gw)$. If $u$ is a positive solution of \eqref{A1} then
	\bel{nontang1} \lim_{x \to y}\myfrac{u(x)}{\BBP^\Gw[\gm](x)}=1 \q non-tangentially, \, \gm-a.e. \ee
Moreover, under the assumptions of theorem A, there holds
	\bel{nontang2} \lim_{x \to y}\myfrac{u(x)}{U_\gm(x)}=1 \q non-tangentially, \, \gm-a.e. \ee
\es
\Proof Put $v_\gm=\BBP^\Gw[\gm]-u$ then $v_\gm>0$ and $-\Gd v_\gm=H \circ u \geq 0$ in $\Gw$. It means $v_\gm$ is a positive superharmonic function in $\Gw$. By Riesz Representation Theorem (see \cite{Ma}), $v_\gm$ can be written under the form $ v_\gm= v_h + v_p $ where $v_h$ is a nonnegative harmonic function and $v_p$ is a potential. Since the boundary trace of $v_\gm$  is a zero measure, it follows that the boundary trace of $v_h$ and $v_p$ is zero measure. Therefore $v_h = 0$ in $\Gw$ and $v_\gm=v_p$. By \cite[Theorem 2.11 and Lemma 2.13]{Ma}, we derive \eqref{nontang1}. Then \eqref{nontang2} follows straightforward from \rcor{max} and \eqref{nontang1}. \qed \medskip

\subsection{Moderate solutions and boundary trace}
In this section we study the notion of boundary trace of positive solutions of \eqref{A0}. We start with some notations.
\bdef{Mtrace} Let $u \in W^{1,s}_{loc}(\Gw)$ for some $s>1$. We say that $u$ possesses an \emph{M-boundary trace} on $\prt \Gw$ if there exists $\gm \in \GTM(\prt \Gw)$ such that, for every uniform $C^2$ exhaustion $\{D_n\}$ (see  \cite[Definition 1.3.1]{MVbook}) and every $\gf \in C(\ovl \Gw)$,
\bel{Mtrace} \lim_{n \to \infty}\myint{\prt D_n}{}u|_{\prt D_n}\gf dS = \myint{\prt \Gw}{}\gf d\gm. \ee
The M-boundary trace of $u$ is denoted by $\tr(u)$.

Let $A$ be a relatively open subset of $\prt \Gw$. A measure $\gm \in \GTM(A)$ is the M-boundary trace of $u$ on $A$ if \eqref{Mtrace} holds for every $\gf \in C(\ovl \Gw)$ such that $\supp \gf \sbs \sbs \Gw \cup A$. In case of positive functions, the definition can be extended to include positive Radon measure on $A$.
\es

Characterization of moderate solutions (see \rdef{moderate}) is given in the next result.

\bth{moderate} Let $u$ be a positive solution of \eqref{A0}. Then the following statements are equivalent:

(i)  $u$ is bounded from above by an harmonic function in $\Gw$.

(ii) $u$ is moderate.

(iii) $u$ possesses an M-boundary trace denoted by $\gm$

(iv) $u$ is a solution of \eqref{A1}.

(v) $u \in L^1(\Gw)$, $H \circ u \in L^1_\rho(\Gw)$ and the integral formulation \eqref{A2} holds where $\mu=\tr(u)$.
\es
\Proof (i) $\Lra$ (ii). Suppose $u \leq U$ where $U$ is a positive harmonic function. By Herglotz's theorem, $U$ admits an M-boundary trace and therefore
$$ \lim_{\gd \to 0}\myint{\Gs_\gd}{}UdS <\infty. $$
It follows that $u \in L^1(\Gw)$ and
$$ \sup_{0<\gd<\gd_0}\myint{\Gs_\gd}{}u dS < \infty. $$
Consequently there exist a sequence $\{\gd_n\}$ converging to zero and a measure $\gm \in \GTM^+(\prt \Gw)$ such that
$$ \lim_{n \to \infty}\myint{\Gs_{\gd_n}}{}u\gf dS = \myint{\prt \Gw}{}\gf d\gm $$
for every nonnegative function $\gf \in C(\ovl \Gw)$. Since $u$ is a solution of \eqref{A0},
\bel{M1} -\myint{D_\gd}{}u \Gd \gz dx + \myint{D_\gd}{}(H \circ u)\gz dx = -\myint{\Gs_\gd}{}u \myfrac{\prt \gz}{\prt {\bf n}}dS \ee
for every $\gz \in C_0^2(\ovl D_\gd)$ and $\gd \in (0,\gd_0)$. Given $\gf \in C^2(\ovl \Gw)$, then there exists a sequence $\{\vgf_n\}$ and a function $\vgf$ such that
\bel{choice} \BA{lll} \vgf_n \in C_0^2(\ovl D_{\gd_n}), \q \myfrac{\prt \vgf_n}{\prt {\bf n}}|_{\Gs_{\gd_n}}=\gf, \q
\vgf \in C_0^2(\ovl \Gw), \q \myfrac{\prt \vgf}{\prt {\bf n}}|_{\prt \Gw} = \gf, \\[4mm]
\norm{\vgf_n}_{C^2(\ovl D_{\gd_n})} < c\norm{\gf}_{C^2(\ovl \Gw)}, \q \vgf_n \leq c \rho_n\gf,\\[4mm]
 \vgf_n/\rho_n \to \vgf/\rho \text{ in } C^2_{loc}(\Gw).
\EA \ee
The constant $c$ is independent of $\gf$ and $n$, but depends on the exhaustion. Consider \eqref{M1} with $\gd=\gd_n$ and $\zeta=\vgf_n$. We see that the first and third terms in \eqref{M1} converge when $n \to \infty$. Therefore the second term converges and we get
\bel{M2} -\myint{\Gw}{}u \Gd \vgf dx + \myint{\Gw}{}(H \circ u)\vgf dx = -\myint{\prt \Gw}{} \myfrac{\prt \vgf}{\prt {\bf n}}d\gm. \ee
By choosing $\gf=1$ in $\Gw$ and $\vgf=\rho$ in $\Gw_{\gd_0/2}$, we deduce that $H \circ u \in L^1_{\rho}(\Gw)$. \medskip

\noindent{(ii) $\Lra$ (iii).} Put $v=u+ \BBG^\Gw[H \circ u]$ then $v$ is a positive harmonic function. Therefore $v$ possesses an M-boundary trace $\gm$. Since $\tr(\BBG^\Gw[H \circ u])=0$, it follows that $\tr(u)=\gm$.  \medskip

\noindent{(iii) $\Lra$ (iv).} The implication is obvious. \medskip

\noindent{(iv) $\Lra$ (v).} Let $\{D_n\}$ be a uniform $C^2$ exhaustion of $\Gw$. For every $n$, denote by $U_n$ the harmonic function in $D_{n}$ such that $U_n=u$ on $\prt D_{n}$. By the comparison principle, $u \leq U_n$ on $D_{n}$. The sequence $\{U_n\}$ converges to a positive harmonic function $U$ which dominates $u$ in $\Gw$. Since $u$ possesses an M-boundary trace $\gm$, it follows that $U$ admits an M-boundary trace $\gm$. Hence $U \in L^1(\Gw)$ and consequently $u \in L^1(\Gw)$. By proceeding as above, we deduce that $H \circ u \in L^1_\rho(\Gw)$. Let $\gf \in C^2(\ovl \Gw)$ and let $\vgf$ and $\{\vgf_n\}$ as in \eqref{choice} with $D_{\gd_n}$ replaced by $D_n$. We have
$$ -\myint{D_n}{}u\Gd \vgf_n dx + \myint{D_n}{}(H \circ u)\vgf_n dx = -\myint{\prt D_n}{}u\gf dS. $$
As $\{\vgf_n/\rho\}$ and $\{\Gd \vgf_n\}$ are bounded sequences converging to $\vgf/\rho$ and $\Gd \vgf$ respectively and $\tr(u)=\gm$, by letting $n \to \infty$, we obtain \eqref{A2}. \medskip

\noindent{(v) $\Lra$ (i).} The implication follows from the estimate $u \leq \BBP^\Gw[\mu]$ in $\Gw$. \qed \medskip

Motivated by the above result, we introduce the following definition.

\bdef{bdtrace}  Let $u$ be a positive solution of \eqref{A0}. A point $y \in \prt \Gw$ is \emph{regular} relative to $u$ if there is a neighborhood $Q$ of $y$ such that
$$ \myint{Q \cap \Gw}{}(H \circ u) \rho\, dx <\infty. $$
Otherwise we say that $y$ is a \emph{singular point} relative to $u$.

The set of regular points is denoted by $\CR(u)$, while the set of singular points is denoted by $\CS(u)$.
\es

\Remark Clearly $\CR(u)$ is relatively open.

The next result can be obtained by combining the argument in the proof of \cite[Theorem 3.1.8]{MVbook} and \rth{moderate}.

\bth{bdtrace}  Let $u$ be a positive solution of \eqref{A0}. Then

(i) $u$ has an M-boundary trace on $\CR(u)$ given by a positive Radon measure $\gm$. Hence
$$ \lim_{\gd \to 0}\myint{\Gs_\gd}{}u\gf dS = \myint{\prt \Gw}{}\gf d\gm $$
for every $\gf \in C(\ovl \Gw)$ such that $\supp \gf \sbs \CR(u)$.

(ii) A point $y \in \prt \Gw$ is singular relative to $u$ if and only if for every $r>0$,
\bel{singpoint} \limsup_{\gd \to 0}\myint{B_r(y) \cap \Gs_\gd}{}udS = \infty. \ee
\es

\Remark From the above results, we see that $u$ is a moderate solution if and only if $\CS(u)=\emptyset$. \medskip

Next we give some results concerning the minimum and the maximum of two positive solutions.
\blemma{subsuper} Let $u_1$ and $u_2$ be two positive solutions of \eqref{A0}. Then $\max(u_1,u_2)$ and $\min(u_1,u_2)$ are respectively a subsolution and a supersolution of \eqref{A0}. Assume in addition that $\tr(u_i)=\gm_i \in \GTM^+(\prt \Gw)$, $i=1,2$. Then $ \tr(\max(u_1,u_2))=\max(\gm_1,\gm_2)$ and $\tr(\min(u_1,u_2))=\min(\gm_1,\gm_2)$.
\es
\Proof Put $v=\max(u_1,u_2)=(u_1-u_2)_+ + u_2$. Since $u_i \in W^{1,s}(\Gw)$ for some $s >1$, it follows that $v \in W^{1,s}(\Gw)$ and
\bel{max1'} \nabla v = \left\{ \BA{ll} \nabla u_1 \q \text{if } u_1 > u_2 \\ \nabla u_2 \q \text{if } u_1 \leq u_2 \EA \right. \q\text{a.e. in } \Gw. \ee
By Kato's inequality (see \cite{MVbook}),
\bel{Kato} \BA{lll} \Gd v &= \Gd (u_1-u_2)_+ + \Gd u_2 \leq sign_+(u_1-u_2)\Gd(u_1-u_2) + \Gd u_2 \\
&= sign_+(u_1-u_2)(H \circ u_1 - H \circ u_2) + H \circ u_2.
\EA \ee
Combining \eqref{max1'} and \eqref{Kato} implies $-\Gd v + H \circ v \leq 0$ in $\Gw$, namely $v$ is a subsolution of \eqref{A0}. Similarly, $\min(u,v)$ is a supersolution of \eqref{A0}.

It follows from \rth{moderate} that $u_i \leq \BBP^\Gw[\gm_i]$, $i=1,2$. Hence
$$v\leq \max(\BBP^\Gw[\gm_1],\BBP^\Gw[\gm_2]).$$
Consequently, $\tr(v)=\max(\gm_1,\gm_2)$. Since
$ \min(u_1,u_2)=u_1+u_2-\max(u_1,u_2) $, it follows that
$$\tr(\min(u_1,u_2))=\gm_1+\gm_2-\max(\gm_1,\gm_2)=\min(\gm_1,\gm_2).$$ \qed

As a consequence, we obtain
\bcor{disjoint} Let $u_i$, $i=1,2$ be positive solutions of \eqref{A0} such that $\tr(u_i)=\gm_i \in \GTM^+(\prt \Gw)$.
Then there exists a minimal solution $\ovl w$ dominating $\max(u_1,u_2)$. This solution  satisfies
\bel{1,2} \max(\mu_1,\mu_2)\leq \tr(\ovl w) \leq \mu_1 +\mu_2.\ee
There exists a nonnegative maximal solution $\unl w$ dominated by $\min(u_1,u_2)$. This solution satisfies
\bel{1+2} \tr(\unl w) \leq \min(u_1,u_2). \ee

 If, in addition, $\supp\gm_1 \cap \supp\gm_2 = \emptyset$ then  $\tr(\ovl w)=\gm_1+\gm_2$. In this case, there exists no positive solution dominated by $\min(u_1,u_2)$.
\es
\Proof For every $\gd \in (0,\gd_0)$, denote by $w:=\ovl w_\gd$ the solution of
\bel{maxdelta} \left\{ \BA{lll} -\Gd w + H \circ w &=0 \qq &\text{in } D_\gd \\
\phantom{-\Gd w + H \circ, } w &= \max(u_1,u_2) &\text{on } \Gs_\gd .
\EA \right. \ee
By the comparison principle, $\max(u_1,u_2) \leq \ovl w_{\gd} \leq \BBP^\Gw[\gm_1+\gm_2]$ in $D_\gd$. Consequently, the sequence $\{\ovl w_{\gd}\}$ is increasing and bounded from above by $\BBP^\Gw[\gm_1+\gm_2]$. Therefore, $\ovl w:=\lim_{\gd \to 0}\ovl w_{\gd}$ is a solution of \eqref{A0} satisfying
$$\max(u_1,u_2) \leq \ovl w \leq \BBP^\Gw[\gm_1+\gm_2]$$
in $\Gw$. By \rth{moderate} $\ovl w$ admits a boundary trace and \eqref{1,2} holds.

If $w$ is a solution of \eqref{A0} dominating $\max\{u_1,u_2\}$ then by the comparison principle, $w \geq \ovl w_\gd$ for every $\gd>0$. It follows that $w \geq \ovl w$ and therefore $\ovl w$ is a minimal solution dominating $\max(u_1,u_2)$.

For every $\gd \in (0,\gd_0)$, denote by $w:=\unl w_\gd$ the solution of
\bel{mindelta} \left\{ \BA{lll} -\Gd w + H \circ w &=0 \qq &\text{in } D_\gd \\
\phantom{-\Gd w + H \circ ,} w &= \min(u_1,u_2) &\text{on } \Gs_\gd .
\EA \right. \ee
The sequence $\{\unl w_\gd\}$ is decreasing and converges, as $\gd \to 0$, to a function $\unl w$ which is a solution of \eqref{A0} such that $0 \leq \unl w \leq \min(u_1,u_2)$ in $\Gw$. As above, one can show that $\unl w$ is the maximal solution dominated by $\min(u_1,u_2)$ and \eqref{1+2} holds.

If $\supp\gm_1 \cap \supp\gm_2 = \emptyset$ then $\tr(\max(u_1,u_2)) = \gm_1+\gm_2$ and $\tr(\min(u_1,u_2))=0$. Therefore $\tr(\ovl w) = \gm_1+\gm_2$ and $\tr(\unl w)=0$. Thus $\unl w \equiv 0$.  \qed

\section{Isolated boundary singularities} If $u$ is a solution of \eqref{A0} in $\Gw$ with an isolated singularity at a point $A\in \bdw$, we shall assume that the set of coordinates is chosen so that $A$ is the origin.
\subsection{Weakly singular solutions} We start with some a-priori estimates regarding solutions with an isolated singularity.
\blemma{estfunct2} Assume $u\in C(\overline \Gw\setminus \{0\})\cap C^2(\Gw)$ is a nonnegative solution of \eqref{A0} in $\Gw$ vanishing on $\prt\Gw\setminus \{0\}$. \medskip

(i) Assume $H$ satisfies \eqref{multi}. Then
\bel{C8} u(x)\leq \Gl_1|x|^{-\gb_1}\qq\forall x\in\Gw, \ee
\bel{C8a} |\nabla u(x)| \leq \Gl_3\abs{x}^{-\gb_1-1} \forevery x \in \Gw, \ee
\bel{C8b}  u(x)| \leq \tl \Gl_3\rho(x)\abs{x}^{-\gb_1-1} \forevery x \in \Gw \ee
where $\Gl_1$ is defined in \eqref{Lambda1}, $\Gl_3=\Gl_3(N,p,q,\Gw)$ and $\tl \Gl_3=\tl \Gl_3(N,p,q,\Gw)$ . \medskip

(ii) Assume $H$ satisfies \eqref{add}. Then
\bel{C8-add} u(x)\leq \Gl_2|x|^{-\gb_2}\qq\forall x\in\Gw, \ee
\bel{C8a-add} |\nabla u(x)| \leq \Gl_4\abs{x}^{-\gb_2-1} \forevery x \in \Gw, \ee
\bel{C8b-add}  u(x)| \leq \tl \Gl_4\rho(x)\abs{x}^{-\gb_2-1} \forevery x \in \Gw \ee
where $\Gl_2=\Gl_2(p,q)$, $\Gl_4=\Gl_4(N,p,q,\Gw)$ and $\tl \Gl_4=\tl \Gl_4(N,p,q,\Gw)$.
\es
\Proof We deal only with the case where $H$ satisfies \eqref{multi} since the case $H$ satisfies \eqref{add} can be treated in a similar way. For $\ge>0$, we set
$$P_\ge(r)=\left\{\BA{ll}
0&\text{if }r\leq\ge\\
\frac{-r^4}{2\ge^3}+\frac{3r^3}{\ge^2}-\frac{6r^2}{\ge}+5r-\frac{3\ge}{2}\q&\text{if }\ge<r<2\ge\\
r-\frac{3\ge}{2}&\text{if } r \geq 2\ge
\EA\right.
$$
and let $u_\ge$ be the extension of $P_\ge(u)$ by zero outside $\Gw$. There exists $R_0$ such that $\Gw\subset B_{R_0}$. Since $0\leq P'_\ge(r)\leq 1$ and $P_\ge$ is convex, $u_\ge \in C^2(\BBR^N \sms\{0\})$ and it satisfies $-\Gd u_\ge+u_\ge^p|\nabla u_\ge|^q\leq 0$. Furthermore $u_\ge$ vanishes in $ B^c_{R_0}$. For $\gd>0$, we set
$$U_{\gd}(x)=\Gl_1(|x|-\gd)^{-\gb_1}\qq\forall x\in \BBR^N\setminus B_\gd,
$$
 then $-\Gd U_{\gd}+U_{\gd}^p|\nabla U_{\gd}|^q\geq 0$ in $B_{R_0} \sms B_\gd$. Since $u_\ge$ vanishes on $\prt B_{R_0}$ and is finite on $\prt B_\gd$, by the comparison principle, $u_\ge\leq U_{\gd}$ in $B_{R_0} \sms {\ovl B_\gd}$. Letting successively $\gd \to 0$ and $\ge\to 0$ yields to \eqref{C8}.

For $\ell>0$, define $T^1_\ell[u](x)=\ell^{\gb_1}u(\ell x)$, $ x \in \Gw^\ell:=\ell^{-1}\Gw$. If $x_0\in\Gw$, we set $r_0=|x_0|$ and $u_{r_0}(x)=T^1_{r_0}[u](x)$. Then $u_{r_0}$ satisfies \eqref{A0} in $\Gw^{r_0}=r_0^{-1}\Gw$. By \eqref{C8},
$$\max\{|u_{r_0}(x)|: (B_\frac{3}{2} \sms B_\frac{1}{2}) \cap \Gw^{r_0} \} \leq 2^{\gb_1}\Gl_1.$$
By regularity results \cite[Theorem 1]{Lib}, there exists $\Gl_3=\Gl_3(N,\Gw,p,q)$ such that $$\max\{|\nabla u_{r_0}(x)|:(B_\frac{5}{4}\sms B_\frac{3}{4}) \cap \Gw^{r_0} \} \leq \Gl_3.$$
In particular, $|\nabla u_{r_0}(x)| \leq \Gl_3$ with $|x|=1$. Hence $|\nabla u(x_0)| \leq \Gl_3|x_0|^{-\gb_1-1}$.

Finally, \eqref{C8b} follows from \eqref{C8} and \eqref{C8a}. \qed \medskip

An uniqueness result for \eqref{A1} can be obtained if $\gm$ is a bounded measure concentrated at a point on $\prt \Gw$. We assume that the point is the origin.
\bth{unique} Assume either $H$ satisfies \eqref{multi} with $0<N(p+q-1)<p+1$ or $H$ satisfies \eqref{add} with $m_{p,q}<p_c$ where $p_c$ and $m_{p,q}$ are given in \eqref{pc} and \eqref{mpq} respectively. Then for every $k >0$, there exists a unique solution,  denoted by $u^\Gw_{k,0}$, of the problem
           \bel{APE0} \left\{ \BA{lll}
	- \Gd u + H \circ u &= 0 \qq &\text{in } \Gw\\
	\phantom{- \Gd u + H \circ,}
	u &= k\gd_0 &\text{on } \prt \Gw.
	\EA \right. \ee
Moreover,
\bel{ab} u^\Gw_{k,0}(x)=k(1+o(1))P^\Gw(x,0) \quad \text{as } x \to 0. \ee
Consequently the mapping $k \mapsto u^\Gw_{k,0}$ is increasing.
\es

The existence of a solution to \eqref{APE0} is guaranteed by Theorem A. The proof of uniqueness is based on the following lemma.

\blemma{APE-grad} Under the assumption of \rth{unique}, let $u$ be a solution to \eqref{APE0}. Then
\bel{AB4} \BBG^{\Gw}[H \circ u](x)= o(P^\Gw(x,0)) \quad \text{as } x \to 0. \ee
\es
\Proof We prove \eqref{AB4} in the case $H$ satisfies \eqref{multi}. The case $H$ satisfies \eqref{add} can be treated in a similar way.

Since $u$ is a solution of \eqref{APE0}, it follows from the maximum principle that
$$u(x) \leq kP^\Gw(x,0) \leq kc_N|x|^{1-N} \forevery x \in \Gw $$
where $c_N$ is a positive constant depending on $N$ and $\Gw$. By adapting argument in the proof of \rlemma{estfunct2}, we obtain
\bel{APE-grad1} \abs{\nabla u(x)} \leq  \Gl_5 k\abs{x}^{-N} \forevery x \in \Gw \ee
where $\Gl_5$ is a positive constant depending on $N,p,q,\Gw$. Consequently,
\bel{G1} \BBG^\Gw[H \circ u](x) \leq  c_8\myint{\Gw}{}G^\Gw(x,y)\abs{y}^{-(N-1)p-Nq}dy \forevery x \in \Gw. \ee
By \cite{MVbook}, there exists $c_9=c_9(N,\Gw)$ such that, for $\vge_0 \in (0,1)$,
	\bel{AB5} \BA{lll} G^\Gw(x,y) \leq c_9\rho(x)\rho(y)^{1-\vge_0}\abs{x-y}^{\vge_0-N} \forevery x,y\in \Gw, x \neq y,
	\EA \ee
This, joint with \eqref{G1}, implies
	\bel{AB6} \BA{lll}\BBG^{\Gw}[H \circ u](x)
\leq c_{10}\abs{x}^NP^\Gw(x,0)\myint{\BBR^N}{}\abs{x-y}^{\vge_0-N}\abs{y}^{1-(N-1)p-Nq-\vge_0}dy. \EA \ee
We fix $\vge_0$ such that $0<\vge_0<\min\left\{1,N+1-(N-1)p-Nq\right\}$. By the following identity (see \cite{LiLo}),
	\bel{AB7} \myint{\BBR^N}{}\abs{x-y}^{\vge_0-N}\abs{y}^{1-(N-1)p-Nq-\vge_0}dy = c_{11}\abs{x}^{N+1-(N-1)p-Nq} \ee
where $c_{11}=c_{11}(N,\vge_0)$, we obtain
$$\BBG^{\Gw}[H \circ u](x) \leq \ga_1c_{10} c_{11}\abs{x}^{N+1-(N-1)p-Nq}P^\Gw(x,0).$$
Since $N+1-(N-1)p-Nq>0$, by letting $x \to 0$, we obtain \eqref{AB4}. \qed \medskip

\noindent{\bf Proof of \rth{unique}.} Let $u_1$ and $u_2$ be two solutions of \eqref{APE0} then $ u_i(x)=k\,P^\Gw(x,0)-\BBG^{\Gw}[H \circ u_i](x)$.  From \eqref{AB4}, we obtain
	\bel{AB9} u_i(x)=k(1+o(1))P^\Gw(x,0) \quad \text{as } x \to 0. \ee
By the comparison principle, we deduce $u_1=u_2$. The monotonicity of $k\mapsto u^\Gw_{k,0}$ follows from \eqref{ab} and the comparison principle. \qed \medskip

Since $\prt \Gw$ is of class $C^2$, there exists $r_0>0$ such that for every $y \in \prt \Gw$, $B_{r_0}(y-r_0\textbf{n}_y) \sbs \Gw$ where $\textbf{n}_y$ is the outward normal unit vector at $y$.
\blemma{Hopf} Assume $H$ satisfies either \eqref{multi} with $p+q>1$ and $q<2$ or $H$ satisfies \eqref{add} with $p>1$, $1<q<2$. Let $u \in C^2(\Gw)$ be a positive solution of \eqref{A0}. Let $y \in \prt \Gw$ be such that $u$ is continuous at $y$ and $u(y)=0$. Then
\bel{liminf} \liminf_{x \to y}\frac{u(x)}{|x-y|}>0 \quad \text{n.t.} \ee
\es
\Proof We only deal with the case $H$ satisfies \eqref{multi} since the case $H$ satisfies \eqref{add} can be treated in a similar way. Put $z=y-r_0\textbf{n}_y$ and set
$$ v(x)=e^{-\ga|x-z|^2}-e^{-\ga r_0^2} \quad x \in B_{r_0}(z) \sms B_{r_0/2}(z)$$
where $\ga>0$ will be determined latter on. Then, in $B_{r_0}(z) \sms B_{r_0/2}(z)$,
$$ \BA{lll} -\Gd v + v^p|\nabla v|^q&=2\ga Ne^{-\ga r^2} - 4\ga^2 r^2 e^{-\ga r^2} + (e^{-\ga r^2}-e^{-\ga r_0^2})^p(2\ga)^q r^q e^{-\ga q r^2} \\
&\leq 2\ga(N-2\ga r^2 + 2^{q-1}\ga^{\ga-1}r^q)e^{-\ga r^2} \\
&\leq \ga(2N-\ga r_0^2 + 2^q\ga^{q-1}r_0^q)e^{-\ga r^2}.
\EA $$
Since $q<2$, one can choose $\ga$ large enough such that the last expression is negative. Consequently, $v$ is a subsolution of \eqref{multi}. As $u$ is positive on $\prt B_{r_0/2}(z)$, one can choose $\vge$ small such that $u > \vge v$ on $\prt B_{r_0/2}(z)$. Obviously $u \geq 0 = \vge v$ on $\prt B_{r_0}(z)$. By the comparison principle, $u \geq \vge v$ in $B_{r_0}(z) \sms B_{r_0/2}(z)$. It follows that
$$ \liminf_{x \to y}\frac{u(x)}{|x-y|} \geq -\vge \frac{\prt v}{\prt \textbf{n}}(y) = 2\vge \ga r_0 e^{-\ga r_0^2} > 0.
$$ \qed
\bprop{compare} Under the assumption of \rth{unique}, for every $k>0$, there exists a positive constant $d_k$  depending on $N$, $p$, $q$, $k$ and $\Gw$ such that
	\bel{compare} d_k P^\Gw(x,0)<u^\Gw_{k,0}(x) < k P^\Gw(x,0) \forevery x \in \Gw. \ee
\es
\Proof The second inequality follows straightforward from the comparison principle. In order to prove the first inequality, put $\CA=\{d>0: d\, P^\Gw(.,0) < u^\Gw_{k,0} \q \text{in } \Gw\}$. Suppose by contradiction that $\CA=\emptyset$. Then for each $n \in \BBN$, there exists a point $x_n \in \Gw$ such that
\bel{contra}  n\, u^\Gw_{k,0}(x_n) \leq  P^\Gw(x_n,0). \ee
We may assume that $\{x_n\}$ converges to a point $x^* \in \ovl \Gw$. We deduce from \eqref{contra} that $x^* \notin \Gw$. Thus $x^* \in \prt \Gw$. By \rth{unique}, $x^* \in \prt \Gw \sms B_{\ge}(0)$ for some $\ge>0$. Denote by $\gs_{x_n}$ the projection of $x_n$ on $\prt \Gw$. It follows from \eqref{contra} that
	$$  \myfrac{u^\Gw_{k,0}(\gs_{x_n}) -u^\Gw_{k,0}(x_n)}{\rho(x_n)} \geq \myfrac{1}{n}\myfrac{P^\Gw(\gs_{x_n},0) - P^\Gw(x_n,0)}{\rho(x_n)}. $$
By letting $n \to \infty$, we obtain
$$ \myfrac{\prt u^\Gw_{k,0}}{\prt \textbf{n}}(0)=0 $$
which contradicts \eqref{liminf}. Thus $\CA \neq \emptyset$. Put $d_k=\max\CA$. By combining \eqref{C8} and boundary Harnack inequality, we deduce that $d_k$ depends on $N$, $p$, $q$, $k$ and $\Gw$.  \qed \medskip

\noindent{\bf Proof of Theorem C.} The proof follows from \rth{unique} and \rprop{compare}. \qed \medskip

The next result gives the existence and uniqueness of the weakly singular solution in the case that $\Gw$ is an unbounded domain.
\bth{unbd} Under the assumption of \rth{unique}, let  either $\Gw=\BBR_+^N:=[x_N>0]$ or $\prt\Gw$ be compact with $0\in\prt\Gw$ ($\Gw$ is possibly unbounded). Then there exists a unique solution $u^\Gw_{k,0}$ to \eqref{APE0}.
\es
\Proof If $\prt\Gw$ is compact, for $n \in \BBN$ large enough, $\prt \Gw \sbs B_n(0)$. We set $\Gw_n=\Gw\cap B_n(0)$ and denote by $u^{\Gw_n}_{k,0}$ the unique solution of
 	\bel {Pc+1}\left\{ \BA{lll} -\Gd u + H \circ u  &= 0 \qq &\text{in } \Gw_n \\
 	\phantom{ -\Gd u + H \circ,}
 	u &= k\gd_0 &\text{on } \prt \Gw_n.
 	\EA \right. \ee
Then
	 \bel {Pc+2}u^{\Gw_n}_{k,0}(x)\leq kP^{\Gw_n}(x,0)\qq\forall x\in \Gw_n
	 \ee	
and
\bel{abn} u_{k,0}^{\Gw_n}(x)=k(1+o(1))P^{\Gw_n}(x,0) \quad \text{as } x \to 0. \ee
Thus $\{u^{\Gw_n}_{k,0}\}$ increase to a function $u^*$ which satisfies
	 	\bel {Pc+3}
	u^*(x)\leq kP^{\Gw}(x,0)\qq\forall x\in \Gw.
	\ee
By regularity theory, $\{u^{\Gw_n}_{k,0}\}_n$ converges in $C^1_{loc}(\ovl \Gw \sms \{0\})$ when $n\to\infty$, and  thus $u^*\in C(\overline\Gw\setminus\{0\})$ is a positive solution of \eqref{A0} in $\Gw$ vanishing on $\prt\Gw\setminus\{0\}$. Estimate \eqref{Pc+3} implies that the boundary trace of $u^*$ is a Dirac measure at $0$, which is in fact $k\gd_0$ due to  \eqref{abn}. Uniqueness is follows from the comparison principle.  \qed \medskip

\subsection{Strongly singular solutions} In this section, we establish existence and uniqueness of strongly singular solutions at a boundary point. We assume that the point is the origin. \medskip

\blemma{singusol} Under the assumption of Theorem C, if $v$ is a positive solution of \eqref{A0} and $y \in \CS(v)$ then $v \geq u^\Gw_{\infty,y}$.
\es
\Proof We can suppose that $y$ is the origin. Since $0 \in \CS(v)$, by  \rth{bdtrace} for every $n\in\BBN_*$,
$$ \limsup_{\gb \to 0}\myint{B_{\frac{1}{n}}(0) \cap \Gs_\gd}{}vdS =\infty. $$
Consequently, there exists a sequence $\{\gd_{n,m}\}_{m\in\BBN}$ tending to zero as $m \to \infty$ such that
$$  \lim_{m\to \infty}\int_{\Gs_{\gd_{n,m}}\cap B_{\frac{1}{n}}(0)}^{}v\,dS=\infty. $$
Then, for any $k>0$, there exists $m_k:=m_{n,k}\in\BBN$ such that
 \bel{E20}
 m\geq m_k\Longrightarrow \myint{\Gs_{\gd_{n,m_k}}\cap B_{\frac{1}{n}}(0)}{}vdS\geq k
 \ee
and $m_{k}\to\infty$ when $n\to \infty$.  In particular there exists $t:=t(n,k)>0$ such that
 \bel{E21}
 \myint{\Gs_{\gd_{n,m_k}}\cap B_{\frac{1}{n}}(0)}{}\inf\{v,t\}dS=k.
 \ee
By the comparison principle $v$ is bounded from below in $D_{\gd_{n,m_k}}$ by the solution $w:=w_{\gd_{n,m_k}}$ of
 \bel{E22}\left\{\BA {ll}
 -\Gd w+ H \circ w =0\qq&\text{in }D_{\gd_{n,m_k}}\\[2mm]
 \phantom{ -\Gd w+ H \circ ,}
w=\inf\{v,t\}\qq&\text{on }\Gs_{\gd_{n,m_k}}.
\EA\right. \ee
When $n\to\infty$, $\inf\{v,t(n,k)\}dS$ converges weakly to $k\gd_0$. By \rcor{var-stab} there exists a subsequence, still denoted by $\{w_{\gd_{n,m_k}}\}_n$,  such that $w_{\gd_{n,m_k}}\to u_{k,0}^\Gw$ when $n\to\infty$ where $u_{k,0}^\Gw$ is the unique solution of \eqref{APE0} and consequently
$v\geq u_{k,0}^\Gw$ in $\Gw$. Therefore $v \geq u^\Gw_{\infty,0}$. \qed \medskip

\noindent{\bf Proof of Theorem D.} By Theorem C and \rlemma{estfunct2}, the sequence $\{u^\Gw_{k,0}\}$ is nondecreasing and bounded from above by either $\tl \Gl_3\rho(x)|x|^{-\gb_1-1}$ if $H$ satisfies \eqref{multi} or $\tl \Gl_4\rho(x)|x|^{-\gb_2-1}$ if $H$ satisfies \eqref{add}. Therefore $\{u^\Gw_{k,0}\}$ converges to a function $u^\Gw_{\infty,0}$. By regularity theory, $u^\Gw_{\infty,0}$ is a solution of \eqref{A0} vanishing on $\prt \Gw \sms \{0\}$. Moreover, since $u^\Gw_{\infty,0} \geq u^\Gw_{k,0}$ for every $k>0$, $\CS(u^\Gw_{\infty,0})=\{0\}$ and  therefore $u^\Gw_{\infty,0} \in \CU_0^\Gw$. If $v \in \CU^\Gw_0$ then by \rlemma{singusol}, $v \geq u^\Gw_{\infty,0}$. Thus $u^\Gw_{\infty,0}$ is the minimal element of $\CU^\Gw_0$.  \qed \medskip

For any $\ell>0$ and any solution of \eqref{A0}, define
\bel{T} \Gw^\ell={\ell^{-1}}\Gw, \q T^1_\ell[u](x)=\ell^{\gb_1}u(\ell x), \q  T^2_\ell[u](x)=\ell^{\gb_2}u(\ell x) \q \forall x \in \Gw^\ell. \ee
\bprop{halfspace} Let $v \in C(\ovl \Gw \sms \{0\}) \cap C^2(\Gw)$ be a nonnegative solution of \eqref{A0} vanishing on $\prt \Gw \sms \{0\}$. \smallskip

\noindent{(1)} Assume $H$ satisfies \eqref{multi}. For each $\ell>0$, put $v_\ell(x)=T^1_\ell[v](x)$. Then, up to a subsequence, $\{v_\ell\}$ converges in $C^1_{loc}(\ovl {\BBR_+^N} \sms \{0\})$, as $\ell \to 0$, to a solution of
\bel{half1} -\Gd u +  u^p |\nabla u|^q = 0 \text{ in } \BBR^N_+, \q u = 0 \text{ on } \prt \BBR^N_+\sms\{0\}. \ee
\noindent{(2)} Assume $H$ satisfies \eqref{add}. For each $\ell$, put $v_\ell(x)=T^2_\ell[v](x)$. \smallskip

(i) If $p=\frac{q}{2-q}$ then, up to a subsequence, $\{v_\ell\}$ converges in $C^1_{loc}(\ovl {\BBR_+^N} \sms \{0\})$, as $\ell \to 0$, to a solution of
\bel{half2} -\Gd u +   u^p + |\nabla u|^q = 0 \text{ in } \BBR^N_+, \q u = 0 \text{ on } \prt \BBR^N_+ \sms\{0\}. \ee

(ii) If $p>\frac{q}{2-q}$ then, up to a subsequence, $\{v_\ell\}$ converges in $C^1_{loc}(\ovl {\BBR_+^N} \sms \{0\})$, as $\ell \to 0$, to a solution of
\bel{half3} -\Gd u +   u^p = 0 \text{ in } \BBR^N_+, \q u = 0 \text{ on } \prt \BBR^N_+\sms\{0\}. \ee

(iii) If $p<\frac{q}{2-q}$ then, up to a subsequence, $\{v_\ell\}$ converges in $C^1_{loc}(\ovl {\BBR_+^N} \sms \{0\})$, as $\ell \to 0$, to a solution of
\bel{half4} -\Gd u +  |\nabla u|^q = 0 \text{ in } \BBR^N_+, \q u = 0 \text{ on } \prt \BBR^N_+\sms\{0\}. \ee
\es
\Proof We first notice that if $H$ satisfies either \eqref{multi} or \eqref{add} with $p=\frac{q}{2-q}$ then $v_\ell$ is a solution of \eqref{A0} in $\Gw^\ell$ which vanishes on $\prt \Gw^\ell \sms \{0\}$. If $H$ satisfies \eqref{add} with $p>\frac{q}{2-q}$ then $v_\ell$ satisfies
\bel{half3l} -\Gd v_\ell +   v_\ell^p + \ell^{\frac{p(2-q)-q}{p-1}}|\nabla v_\ell|^q= 0 \text{ in } \Gw^\ell, \q v_\ell = 0 \text{ on } \prt \Gw^\ell \sms \{0\}. \ee
If $H$ satisfies \eqref{add} with $p<\frac{q}{2-q}$ then $v_\ell$ satisfies
\bel{half4l} -\Gd v_\ell +  \ell^{\frac{q-(2-q)p}{q-1}}v_\ell^p + |\nabla v_\ell|^q= 0 \text{ in } \Gw^\ell, \q v_\ell = 0 \text{ on } \prt \Gw^\ell \sms \{0\}. \ee

Next, it follows from \rlemma{estfunct2} and \cite[Theorem 1]{Lib} that for every $R>1$ there exist positive numbers $M=M(N,p,q,R)$ and $\gg=\gg(N,p,q) \in (0,1)$ such that
\bel{lie} \BA{ll} \sup\{\abs{v_\ell(x)}+\abs{\nabla v_\ell(x)}: x \in \Gg_{R^{-1},R} \cap \Gw^{\ell}\} \\[3mm]
   \phantom{qqqqq}
   + \sup\left\{\myfrac{\abs{\nabla v_\ell(x) - \nabla v_\ell (y)}}{\abs{x-y}^\gg}: x,y \in \Gg_{R^{-1},R} \cap \Gw^{\ell} \right\} \leq M \EA \ee
where $\Gg_{t_1,t_2}:=B_{t_2}(0)\sms B_{t_1}(0)$ with $0<t_1<t_2$. Thus there exists a sequence $\{\ell_n\}$ and a function $v^* \in C^1(\ovl{\BBR^N_+} \sms \{0\})$ such that $\{v_{\ell_n}\}$ converges to $v^*$ in $C^1_{loc}(\ovl{\BBR^N_+} \sms \{0\})$ which is a solution of
    \bel{halfspace} \left\{ \BA{lll} -\Gd v +  H \circ v & = 0 \qq &\text{in } \BBR^N_+ \\
    \phantom{-\Gd v +  H \circ ,}
    v &= 0 &\text{in } \prt \BBR^N_+ \sms \{0\}\EA \right. \ee
Moreover,
    \bel{t2} \lim_{n \to \infty}(\sup\{|(v_{\ell_n}-v^*)(x)|+|\nabla (v_{\ell_n}-v^*)(x)|: x \in \Gg_{R^{-1},R} \cap \Gw^{\ell_n}\})=0. \ee \qed
\bprop{halfspace-c} Let $v=u^\Gw_{\infty,0}$ and $\{v_\ell\}$ be defined as in \rprop{halfspace}. Then, up to a subsequence,  $\{v_\ell\}$ converges to a strongly singular solution of
\bel{v-k}  \left\{ \BA{lll} (\ref{half1}) \text{ if } H \text{ satisfies } (\ref{multi}) \\
(\ref{half2}) \text{ if } H \text{ satisfies } (\ref{add}) \text{ with } p=\frac{q}{2-q} \\
(\ref{half3}) \text{ if } H \text{ satisfies } (\ref{add}) \text{ with } p>\frac{q}{2-q} \\
(\ref{half4}) \text{ if } H \text{ satisfies } (\ref{add}) \text{ with } p<\frac{q}{2-q}
\EA \right. \ee
Moreover, if $H$ satisfies \eqref{add} and $p \neq \frac{q}{2-q}$ then the whole sequence $\{v_\ell\}$ converges to the strongly singular solution of \eqref{half3} if $p>\frac{q}{2-1}$ or it converges to the strongly singular solution of \eqref{half4} if $p<\frac{q}{2-1}$.
\es
\Proof Since $v_\ell \geq u^{\Gw_\ell}_{k,0}$ for every $\ell>0$ and $k>0$, $v^* \geq u^{\BBR^N_+}_{k,0}$ for every $k>0$ ($v^*$ is given in the proof of \rprop{halfspace}). Therefore $v^*$ is a strongly singular solution.

If $H$ satisfies \eqref{add} with $p \neq \frac{q}{2-q}$ then by the uniqueness of the strongly singular solution of \eqref{half3} and \eqref{half4} (see \cite{MV4} and \cite{NV}), we get the conclusion. \qed \medskip

Denote by $\CE_i$ $(i=1,\ldots,4$) the set of positive solutions in $C^2(S^{N-1}_+)$ of
\bel{PH1} -\Gd' \gw +  F_i(\gw,\nabla' \gw)=0 \text{ in } S_+^{N-1}, \q \gw = 0 \text{ on } \prt S_+^{N-1} \ee
where $F_i$ is as in \eqref{i=1...4}.

We next study the structure of $\CE_i$.

\bth{exis-uniq-SS}  (i)  Assume either $H$ satisfies \eqref{multi} with $N(p+q-1) \geq p+1$ or $H$ satisfies \eqref{add} with $m_{p,q} \geq p_c$. Then $\CE_i = \emptyset$  where
\bel{ii}  i=\left\{ \BA{lll} 1 \text{ if } H \text{ satisfies } (\ref{multi}) \\
2 \text{ if } H \text{ satisfies } (\ref{add}) \text{ with } p=\frac{q}{2-q} \\
3 \text{ if } H \text{ satisfies } (\ref{add}) \text{ with } p>\frac{q}{2-q} \\
4 \text{ if } H \text{ satisfies } (\ref{add}) \text{ with } p<\frac{q}{2-q}.
\EA \right. \ee

(ii) Assume either $H$ satisfies \eqref{multi} with $0<N(p+q-1)<p+1$ or $H$ satisfies \eqref{add} with $m_{p,q}<p_c$. Then $\CE_i \neq \emptyset$ with $i$ as in \eqref{ii}.
\es
\Proof Notice that if $H$ satisfies \eqref{add} with $p \neq \frac{q}{2-q}$,  statements (i) and (ii) have been proved in \cite{MV4} and \cite{NV}. More precisely, if $m_{p,q}<p_c$ and $p>\frac{q}{2-q}$ then there exists a unique element $\gw_3$ of $\CE_3$, while if $m_{p,q}<p_c$ and $p<\frac{q}{2-q}$ then there exists a unique element $\gw_4$ of $\CE_4$.

So we are left with the case when $H$ satisfies either \eqref{multi} or $H$ satisfies \eqref{add} with $p=\frac{q}{2-q}$ and we only give the proof for the case $H$ satisfies \eqref{multi}.

(i) Denote by $\vgf_1$ the first eigenfunction of $-\Gd'$ in $W^{1,2}_0(S^{N-1}_+)$, normalized such that $\max_{S^{N-1}_+}\vgf_1=1$, with corresponding eigenvalue $\gl_1=N-1$. Multiplying \eqref{PH1} by $\vgf_1$ and integrating over $S_+^{N-1}$, we get
	$$ \BA{ll}
	\left[N-1-\gb_1(\gb_1+2-N)\right]\myint{S_+^{N-1}}{}\gw\, \vgf_1 dS(\gs)
+ \myint{S_+^{N-1}}{}\gw^p(\gb_1^2\, \gw^2+\abs{\nabla' \gw}^2)^{\frac{q}{2}}\vgf_1 dS(\gs) = 0. \EA $$
{\it Therefore  if $N-1 \geq \gb_1(\gb_1+2-N)$, namely $N(p+q-1)\geq p+1$, then there exists no positive solution of \eqref{PH1}}. \medskip

(ii) The proof is based on the construction of a subsolution and a supersolution to \eqref{PH1}. By a simple computation, we can prove that $\unl \gw:=\gth_1 \vgf_1^{\gth_2}$ is a positive subsolution of \eqref{PH1} with $\gth_1>0$ small and $ 1< \gth_2< \frac{\gb_1\,(\gb_1+2-N)}{N-1}$. Next, it is easy to see that $\ovl \gw=\gth_3$, with $\gth_3>0$ large enough, is a supersolution of \eqref{PH1} and $\ovl \gw>\unl \gw$ in $\ovl S^{N-1}_+$. Therefore by \cite{KaKr} there exists a solution $\gw_1 \in W^{2,m}(S^{N-1}_+)$ (for any $m>N$) to \eqref{PH1} such that $0< \unl \gw \leq \gw_1 \leq \ovl \gw$ in $S_+^{N-1}$. By regularity theory,  $\gw_1 \in C^{2}(\overline{S^{N-1}_+})$.

Similarly, we can show that if $m_{p,q}<p_c$ and $p=\frac{q}{2-q}$ then there exists a function $\gw_2 \in \CE_2$. \qed \medskip

We next show that $\gw_i$ ($i=1,2 $) is the unique element of $\CE_i$.

\bth{unique-sphere}  (i) If $H$ satisfies \eqref{multi} with $N(p+q-1)<p+1$ and $p \geq 1$ then $\CE_1=\{\gw_1\}$.

(ii) If $H$ satisfies \eqref{add} with $m_{p,q}<p_c$ and $p=\frac{q}{2-q}$ then $\CE_2=\{\gw_2\}$.
\es
\Proof We give below only the proof of statement (i); the statement (ii) can be treated in a similar way. Suppose that $\gw_1$ and $\gw_1'$ are two positive different solutions of \eqref{PH1}. Up to exchanging the role of $\gw_1$ and $\gw_1'$, we may assume $ \max_{S^{N-1}_+}\gw_1' \geq \max_{S^{N-1}_+}\gw_1$ and
	$$ \gt_0:=\inf\{\gt>1:\gt\gw_1 > \gw_1' \text{ in } S^{N-1}_+\} >  1. $$
Set $\gw_{1,\gt_0}:=\gt_0 \gw_1$, then $\gw_{1,\gt_0}$ is a positive supersolution to problem \eqref{PH1}. Put $\tl \gw=\gw_{1,\tau_0}-\gw_1' \geq 0$. If there exists $\gs_0 \in S^{N-1}_+$ such that $\gw_{1,\gt_0}(\gs_0)=\gw_1'(\gs_0)>0$ and $\nabla ' \gw_{1,\gt_0}(\gs_0)=\nabla ' \gw_1'(\gs_0)$ then $\tl \gw(\gs_0)=0$ and $\nabla' \tl \gw(\gs_0)=0$. This contradicts the strong maximum principle (see \cite{GT}). If $\gw_{1,\gt_0}>\gw_1'$ in $S^{N-1}_+$ and  there exists $\gs_0 \in \prt S^{N-1}_+$ such that $\frac{\prt \gw_{1,\gt_0}}{\prt \textbf{n}}(\gs_0)=\frac{\prt \gw_1'}{\prt \textbf{n}}(\gs_0)$ then $\tl \gw > 0$ and $\frac{\prt \tl \gw}{\prt \textbf{n}}(\gs_0)=0$. This contradict the Hopf lemma (see \cite{GT}). \qed \medskip

Let $T^1_\ell$ and $T^2_\ell$ be as in \eqref{T}. If $u$ is a solution of \eqref{A0} in $\Gw$ with $H$ as in \eqref{multi} (resp. $H$ as in \eqref{add} and $p=\frac{q}{2-q}$) then $T^1_\ell[u]$ (resp. $T^2_\ell[u]$) is a solution of \eqref{A0} in $\Gw^\ell=\ell^{-1}\Gw$. If $\Gw=\Gw^\ell$ and $u=T^1_\ell[u]$ (resp. $u=T^2_\ell[u]$) for every $\ell>0$ we say that $u$ is a \emph{self-similar solution}.

For $x \in \BBR^N \sms \{0\}$, put $r=|x|$ and $\gs=\frac{x}{r}$.

\bprop{sing}

(i) If $H$ satisfies \eqref{multi} with $N(p+q-1)<p+1$ and $p \geq 1$ then
\bel{uniq7}
r^{\gb_1}u^\Gw_{\infty,0}(x) \to \gw_1(\gs) \qq \text{as } r \to 0, \; x \in \Gw, \; \gs \in S^{N-1}_+
 \ee
 locally uniformly on $S^{N-1}_+$. \smallskip

(ii) If $H$ satisfies \eqref{add} with $m_{p,q}<p_c$ then
\bel{uniq7'}
r^{\gb_2}u^\Gw_{\infty,0}(x) \to \gw_i(\gs) \qq \text{as } r \to 0, \; x \in \Gw, \; \gs \in S^{N-1}_+
 \ee
locally uniformly on $S^{N-1}_+$ where $i\geq2$ is as in \eqref{ii}.
\es
\Proof {\bf Case 1:} $H$ satisfies \eqref{multi}. Since the proof is close to the one of \cite[Proposition 3.22]{NV}, we present only the main ideas.

We first  note that $T^1_\ell[u^{\BBR^N_+}_{\infty,0}]=u^{\BBR^N_+}_{\infty,0}$ for every $\ell>0$. Hence $u^{\BBR^N_+}_{\infty,0}$ is self-similar and satisfies \eqref{uniq7} with $\Gw$ replaced by $\BBR^N_+$.

Next,  let $B$ and $B'$ are two open balls tangent to $\prt \Gw$ at $0$ such that $B \sbs \Gw \sbs G:=(B')^c$. Then
 \bel{uniq12*} u_{\infty,0}^{B^{\ell'}}\leq u_{\infty,0}^{B^{\ell}} \leq u_{\infty,0}^{\BBR^N_+}\leq  u_{\infty,0}^{G^\ell}
 \leq  u_{\infty,0}^{G^{\ell ''}} \qq\forall \,0<\ell\leq\ell',\ell''\leq 1.
\ee
Notice that $u_{\infty,0}^{B^{\ell}}\uparrow \underline u_{\infty,0}^{\BBR^N_+}$ and $u_{\infty,0}^{G^\ell}\downarrow \overline u_{\infty,0}^{\BBR^N_+}$ when $\ell\to 0$ where $ \underline u_{\infty,0}^{\BBR^N_+}$ and $\overline u_{\infty,0}^{\BBR^N_+}$ are positive solutions of \eqref{A0} in $\BBR^N_+$, continuous in $\overline{\BBR^N_+}\setminus\{0\}$ and vanishing on $\prt\BBR^N_+\setminus\{0\}$. By letting $\ell \to 0$ in \eqref{uniq12*}, we obtain
  \bel{uniq13a}
 u_{\infty,0}^{B^{\ell}} \leq \underline u_{\infty,0}^{\BBR^N_+}\leq u_{\infty,0}^{\BBR^N_+}\leq \overline u_{\infty,0}^{\BBR^N_+}\leq  u_{\infty,0}^{G^\ell}\qq\forall\,0<\ell\leq 1.
\ee
Furthermore there also holds for $\ell,\ell'>0$,
  \bel{uniq14}
T^1_{\ell'\ell}[u_{\infty,0}^{B}]=T^1_{\ell'}[T^1_{\ell}[u_{\infty,0}^{B}]]=u_{\infty,0}^{B^{\ell\ell'}} \text{ and }
T^1_{\ell'\ell}[u_{\infty,0}^{G}]=T^1_{\ell'}[T^1_{\ell}[u_{\infty,0}^{G}]]=u_{\infty,0}^{G^{\ell\ell'}}.
\ee
Letting $\ell\to 0$ in \eqref{uniq14} yields
  \bel{uniq15}
\underline u_{\infty,0}^{\BBR^N_+}=T^1_{\ell'}[\underline u_{\infty,0}^{\BBR^N_+}] \text{ and }
\overline u_{\infty,0}^{\BBR^N_+}=T^1_{\ell'}[\overline u_{\infty,0}^{\BBR^N_+}].
\ee
Thus $\underline u_{\infty,0}^{\BBR^N_+}$ and $\overline u_{\infty,0}^{\BBR^N_+}$ are self-similar solutions of \eqref{A0} in $\BBR^N_+$ vanishing on $\prt \BBR_+^N\sms\{0\}$ and continuous in $\overline{\BBR^N_+}\setminus\{0\}$. Therefore they coincide with $u_{\infty,0}^{\BBR^N_+}$.

Finally, since
 \bel{uniq16} u_{\infty,0}^{B^{\ell}}\leq T^1_{\ell}[u_{\infty,0}^{\Gw}]
 \leq  u_{\infty,0}^{G^\ell}\qq\forall \,0<\ell\leq 1,
\ee
by letting $\ell \to 0$ we obtain \eqref{uniq7}. \medskip

\noindent {\bf Case 2:} $H$ satisfies \eqref{add} with $p=\frac{q}{2-q}$. The proof is similar to the one in case 1.

\noindent {\bf Case 3:} $H$ satisfies \eqref{add} with $p>\frac{q}{2-q}$. For any $k>0$ and $\ell>0$, $T^2_\ell[u^\Gw_{k,0}]$ is a solution of  \eqref{half3l} with boundary trace $\ell^{\gb_2+1-N}k\gd_0$. For $k >0$, denote by $Y_k^\Gw$ the unique positive solution of
\bel{Y} - \Gd Y + Y^p = 0 \q \text{in } \Gw \ee
with boundary trace $k\gd_0$ and by $Y_\infty^\Gw$ the unique solution of \eqref{Y} with strong singularity at the origin.

Since $0<\ell<1$ and $p>\frac{q}{2-q}$, by the comparison principle, we get
$$u^{\Gw^\ell}_{\ell^{\gb_2+1-N} k,0} \leq T^2_\ell[u^\Gw_{k,0}] \leq T_\ell^2[Y_k^\Gw] = Y_{\ell^{\gb_2+1-N}k}^{\Gw^\ell} \q \text{ in } \Gw^\ell.$$ By letting $k \to \infty$, we obtain
$$u^{\Gw^\ell}_{\infty,0} \leq T^2_\ell[u^\Gw_{\infty,0}] \leq Y^{\Gw^\ell}_{\infty} \q \text{ in } \Gw^\ell. $$
By \rprop{halfspace-c} and \cite{MV4}, letting $\ell \to 0$ we deduce that
$$ \lim_{\ell \to 0}\ell^{\gb_2}u^\Gw_{\infty,0}(\ell x)=Y^{\BBR^N_+}_{\infty}(x)=|x|^{\gb_2}\gw_3(x/|x|) $$
which implies \eqref{uniq7'} with $i=3$. \medskip

\noindent {\bf Case 4:} $H$ satisfies \eqref{add} with $p<\frac{q}{2-q}$. Denote by $Z_\infty^{\BBR^N_+}$ the positive solution of
\bel{Z} - \Gd Z + |\nabla Z|^q = 0 \q \text{in } \BBR^N_+ \ee
with strong singular at the origin. By proceeding as in case 3 and results in \cite{NV}, we derive
$$ \lim_{\ell \to 0}\ell^{\gb_2}u^\Gw_{\infty,0}(\ell x)=Z^{\BBR^N_+}_{\infty}(x)=|x|^{\gb_2}\gw_4(x/|x|). $$
Thus \eqref{uniq7'} with $i=4$ follows. \qed \medskip

We next construct the maximal strongly singular solution.

\bprop{max}  (i) Assume either $H$ satisfies \eqref{multi} with $0<N(p+q-1)<p+1$ then there exists a maximal element $U_{\infty,0}^{\Gw}$ of $\CU^\Gw_0$. In addition, if $p \geq 1$ then
\bel{max1}
r^{\gb_1}U^\Gw_{\infty,0}(x) \to \gw_1(\gs) \qq \text{as } r \to 0, \; x \in \Gw, \; \gs \in S^{N-1}_+
 \ee
 locally uniformly on $S^{N-1}_+$. \medskip

(ii) If $H$ satisfies \eqref{add} with $m_{p,q}<p_c$ then there exists a maximal element $U_{\infty,0}^{\Gw}$ of $\CU^\Gw_0$ and
\bel{max1-add}
r^{\gb_2}U^\Gw_{\infty,0}(x) \to \gw_i(\gs) \qq \text{as } r \to 0, \; x \in \Gw, \; \gs \in S^{N-1}_+
 \ee
 locally uniformly on $S^{N-1}_+$ where $i \geq2$ is as in \eqref{ii}.
\es
\Proof {\bf Case 1:}  $H$ satisfies \eqref{multi}. \smallskip

\noindent \textit{Step 1: Construction maximal solution}. Since $0<N(p+q-1)<p+1$, there exists a radial solution of \eqref{A0} in $\BBR^N\setminus\{0\}$ of the form
\bel{max2}
U^\dag_1(x)=\Gl^\dag_1 \,|x|^{-\gb_1}\q\text{with }\;\Gl^\dag_1=\left(\myfrac{\gb_1+2-N}{\gb_1^{q-1}}\right)^{\frac{1}{p+q-1}}.
\ee
Therefore,
 $U_1^{*}(x)=\Gl^*_1|x|^{-\gb_1}$ with $\Gl^*_1= \max\{\Gl^\dag_1,\Gl_1\}$ (here $\Gl_1$ is the constant in \eqref{C8}) is a supersolution of \eqref{A0} in $\BBR^N\setminus\{0\}$ and dominates $u$ in $\Gw$. Let $\{\psi_{\ge,n}\}$ with $0<\ge<\max\{|z|:z\in\Gw\}$ be a decreasing smooth sequence on $(\prt \Gw \sms B_\ge(0)) \cup (\Gw \cap \prt B_\ge(0))$ such that
$$ 0\leq \psi_{\ge,n} \leq \Gl^*_1\ge^{-\gb_1}, \quad  \psi_{\ge,n}(x)=\Gl^*_1\ge^{-\gb_1}\quad \text{ if } x \in \Gw \cap \prt B_\ge(0) $$
$$ \psi_{\ge,n}(x)=0 \quad \text{if } x \in \prt \Gw \sms B_\ge(0) \text{ and } \dist(x,\prt B_\ge(0))>\frac{1}{n}. $$
Let $u^\Gw_{\ge,n}$ be the solution of
 \bel{u aprox}
\left\{\BAL
-\Gd u+ H \circ u&=0\q\text {in }\Gw\setminus B_\ge(0)\\
u&=\psi_{\ge,n} \q\text {on }(\prt\Gw\setminus B_\ge(0)) \cup (\Gw \cap \prt B_\ge(0))
\EAL\right.\ee
By the comparison principle, $u^\Gw_{\ge,n} \leq U^{*}_1 $ in $\Gw\sms B_\ge(0)$. Owing to \rcor{var-stab}, $\{u^\Gw_{\ge,n}\}$ converges to the solution $u^\Gw_\ge$ of
  \bel{max4} \left\{\BAL
-\Gd u_\ge + H \circ u &=0 \q \text {in }\Gw\setminus B_\ge(0) \\
u_\ge&=0\q\text {on }\prt\Gw\setminus B_\ge(0)\\\phantom{------,,}
u_\ge&=\Gl^*_1\ge^{-\gb_1}\q\text {on }\Gw\cap \prt B_\ge(0).
\EAL\right.\ee
Consequently, $u^\Gw_{\ge} \leq U^{*}_1$. If  $\ge'<\ge$, for $n$ large enough, $u^\Gw_{\ge',n} \leq u^\Gw_{\ge,n}$, therefore
  \bel{max5}
  u^\Gw_{\ge'}\leq u^\Gw_\ge\leq U^{*}_1\qq\text{in }\;\Gw.
  \ee
Letting $\ge$ to zero, $\{u^\Gw_\ge\}$ decreases and converges to some $U_{\infty,0}^{\Gw}$ which vanishes on  $\prt\Gw\setminus \{0\}$. Therefore $U_{\infty,0}^{\Gw}\in \CU^\Gw_0$. Moreover, there holds
    \bel{max6}
 u_{\infty,0}^{\Gw} \leq  U_{\infty,0}^{\Gw}\leq U^{*}_1(x).
  \ee
If $u \in \CU^\Gw_0$ then $u \leq u^\Gw_{\ge,n}$. Consequently, $u \leq U_{\infty,0}^{\Gw}$. Therefore $ U_{\infty,0}^{\Gw}$ is the maximal element of $\CU^\Gw_0$. \medskip

  \noindent \textit{Step 2: Proof of \eqref{max1}.} Assume $p \geq 1$. From the fact that
    \bel{max7} T^1_\ell[U^{*}_1]=U^{*}_1 \qq\forall\,\ell>0,
    \ee
and \rth{unique-sphere}, we deduce $U_{\infty,0}^{\BBR^N_+} \equiv u_{\infty,0}^{\BBR^N_+}$.

Next, let $B$ and $B'$ are two open balls tangent to $\prt \Gw$ at $0$ such that $B \sbs \Gw \sbs G:=(B')^c$. Note that $$T^1_\ell[u^B_\ge]=u^{B^\ell}_{\frac{\ge}{\ell}} \q \text{and} \q T^1_\ell[u^G_\ge]=u^{G^\ell}_{\frac{\ge}{\ell}} \forevery \ell,\ge>0 $$  where $u^G_\ge$ is the solution of \eqref{max4} in $G \setminus B_\ge(0)$. By letting $\ge \to 0$ we deduce that
\bel{max9}
 T^1_\ell[U_{\infty,0}^{B}]=U_{\infty,0}^{B^{\ell}} \text{ and }\;T^1_\ell[U_{\infty,0}^{G}]=U_{\infty,0}^{G^\ell}.
\ee
Notice that
  \bel{uniq12} U_{\infty,0}^{B^{\ell'}}\leq U_{\infty,0}^{B^{\ell}} \leq U_{\infty,0}^{\BBR^N_+}\leq  U_{\infty,0}^{G^\ell}
 \leq  U_{\infty,0}^{G^{\ell''}}\qq\forall \,0<\ell\leq\ell',\ell''\leq 1
\ee
and
  \bel{uniq13} U_{\infty,0}^{B^{\ell'}}\leq U_{\infty,0}^{B^{\ell}} \leq T^1_\ell[U_{\infty,0}^{\Gw}] \leq  U_{\infty,0}^{G^\ell}
 \leq  U_{\infty,0}^{G^{\ell''}}\qq\forall \,0<\ell\leq\ell',\ell''\leq 1.
\ee
Hence $U_{\infty,0}^{B^{\ell}}\uparrow \underline U_{\infty,0}^{\BBR^N_+}\leq U_{\infty,0}^{\BBR^N_+}$ and
$U_{\infty,0}^{G^\ell} \downarrow \overline U_{\infty,0}^{\BBR^N_+}\geq U_{\infty,0}^{\BBR^N_+}$ as $\ell \to 0$ where  $\underline U_{\infty,0}^{\BBR^N_+}$ and $ \overline U_{\infty,0}^{\BBR^N_+}$ are positive solutions of \eqref{A0} in $\BBR^N_+$ which vanish on $\prt\BBR^N_+\setminus\{0\}$ and endow the same scaling invariance under $T^1_\ell$. Therefore they coincide with $u_{\infty,0}^{\BBR^N_+}$. Letting $\ell \to 0$ in \eqref{uniq13} implies \eqref{max1}.  \medskip

\noindent {\bf Case 2:} $H$ satisfies \eqref{add} with $p=\frac{q}{2-q}$. Since in this case, \eqref{A0} admits a similarity transformation $T_\ell^2$, the proof is similar to the one in case 1. \medskip

\noindent {\bf Case 3:} $H$ satisfies \eqref{add} with $p>\frac{q}{2-q}$. In this case, \eqref{A0} admits no similarity transformation and there is no radial solution of \eqref{A0} in $\BBR^N \sms \{0\}$. We can instead employ a radial supersolution of the form
\bel{max2-add}
U^*_3(x)=\Gl^*_3\,|x|^{-\gb_2}\q\text{with }\;\Gl^*_3={\gb_2(\gb_2+2-N)}^{\frac{1}{p-1}}
\ee
and then we proceed to construct the maximal solution as in case 1. For $\ge>0$, let $u^\Gw_\ge$ be the solution of \eqref{max4}.  Since $u^{\Gw^\ell}_{\frac{\ge}{\ell}} \leq T^2_\ell[u^\Gw_\ge]$, by letting $\ge \to 0$ we obtain $U^{\Gw^\ell}_{\infty,0} \leq T^2_\ell[U^\Gw_{\infty,0}]$. It follows that
$$ u^{\Gw^\ell}_{\infty,0} \leq  U^{\Gw^\ell}_{\infty,0}  \leq  T^2_\ell[U^\Gw_{\infty,0}] \leq  T^2_\ell[Y^\Gw_\infty] =  Y^{\Gw^\ell}_{\infty} $$
where $Y^{\Gw}_{\infty}$ is the unique strongly singular solution of \eqref{Y}. Due to \rprop{halfspace-c} and the uniqueness, we deduce
	$$ \lim_{\ell \to 0}T^2_\ell[U^\Gw_{\infty,0}]=Y^{\BBR^N_+}_{\infty},$$
which, together with the fact $Y^{\BBR^N_+}_{\infty}(x)=|x|^{-\gb_2}\gw_3(x/|x|)$, implies \eqref{max1-add}. \medskip

\noindent {\bf Case 4:} $H$ satisfies \eqref{add} with $p<\frac{q}{2-q}$. The proof is similar to the one in case 3. \qed \medskip

 \rprop{sing} and \rprop{max} show that the minimal solution $u^\Gw_{\infty,0}$ behaves like the maximal solution $U^\Gw_{\infty,0}$ near the origin, which enables us to prove the following result.

\bth {UNI} Assume either $H$ satisfies \eqref{multi} with $N(p+q-1)<p+1$ and $p \geq 1$ or $H$ satisfies \eqref{add} with $m_{p,q}<p_c$. Then $U^\Gw_{\infty,0}=u^\Gw_{\infty,0}$.
  \es
\Proof {\bf Case 1:} $H$ satisfies \eqref{multi} with $p\geq 1$. \smallskip

We represent $\prt\Gw$ near $0$ as the graph of a $C^2$ function $\gf$
defined in $\BBR^{N-1}\cap B_R$ and such that $\gf(0)=0$, $\nabla\gf(0)=0$ and
$$\prt\Gw\cap B_R=\{x=(x',x_N):x'\in \BBR^{N-1}\cap B_R,x_N=\gf(x')\}.$$
We introduce the new variable $y=\Gf(x)$ with $y'=x'$ and $y_N=x_N-\gf(x')$, with corresponding spherical coordinates in $\BBR^N$, $(r,\gs)=(|y|,\frac{y}{|y|})$.

Let $u$ is a positive solution of \eqref{A0} in $\Gw$ vanishing on $\prt\Gw \sms \{0\}$. We set $u(x)=r^{-\gb_1}v(t,\gs)$ with  $t=-\ln r\geq 0$, then a technical computation shows that $v$ satisfies with ${\bf n}=\frac{y}{|y|}$
\begin{equation}\label{uni2}\BA {l}
\left(1+\ge^1_1\right)v_{tt}+\left(2\gb_1+2-N+\ge^1_2\right)v_{t}
+(\gb_1\left(\gb_1+2-N)+\ge^1_3\right)v+\Gd'v\\[3mm]
\phantom{--}
+\langle\nabla'v,\overrightarrow {\ge^1_4}\rangle+\langle\nabla'v_t,\overrightarrow {\ge^1_5}\rangle+
\langle\nabla'\langle \nabla' v,{\bf e}_N\rangle,\overrightarrow {\ge^1_6}\rangle\\[3mm]\phantom{--}
-v^p\abs{(-\gb_1\,v+v_t){\bf n}+\nabla' v+\langle(-\gb_1\,v+v_t){\bf n}+\nabla' v,{\bf e}_N\rangle\overrightarrow\ge^1_7}^q=0,
\EA\end{equation}
on $Q_R:=[-\ln R, \infty)\ti S^{N-1}_{+}$ where $\ge^1_j$ have the following properties
\begin{itemize}
\item $\ge^1_j$ are uniformly continuous functions of $t$ and $\gs\in S^{N-1}$ for $j=1,...,7$,
\item $\ge^1_j$ are $C^1$ functions for $j=1,5,6,7$,
\item $|\ge^1_j(t,.)|\leq c_{12}e^{-t}$  for $j=1,...,7$ and $|\ge^1_{j\,t}(t,.)|+|\nabla'\ge^1_j|\leq c_{12}e^{-t}$  for $j=1,5,6,7$.
\end{itemize}
Moreover $v$ vanishes on $[-\ln R,\infty)\ti \prt S^{N-1}_{+}$.
By \cite[Theorem 4.7]{GV}, there exist constants $c_{13}>0$ and $T>\ln R$ such that
\begin{equation}\label{uni4}\BA {l}
\norm {v(t,.)}_{C^{2,\gg}(\overline {S^{N-1}_{+}})}+\norm {v_t(t,.)}_{C^{1,\gg}(\overline {S^{N-1}_{+}})}
+\norm {v_{tt}(t,.)}_{C^{0,\gg}(\overline {S^{N-1}_{+}})}\leq c_{13}
\EA\end{equation}
for $\gamma\in (0,1)$ and $t \geq T+1$. Moreover
$$ \lim_{t \to \infty}^{}\int_{S^{N-1}_+}^{}(v_t^2 + v_{tt}^2 + |\nabla ' v_t|^2)dS(\gs) = 0. $$
Consequently, the $\gw$-limit set of $v$
$$ \Gg^+(v)=\cap_{\tau \geq 0}\ovl{\cup_{t \geq \tau}v(t,.)}^{C^2(S^{N-1}_+)} $$
is a non-empty, connected and compact subset of the set of $\CE_1$. By the uniqueness of \eqref{PH1}, $\Gg^+(v)=\CE_1=\{\gw_1\}$. Hence $\lim_{t\to \infty}v(t,.)=\gw_1$ in $C^2(\overline {S^{N-1}_{+}})$.\smallskip

By taking $u=u^\Gw_{\infty,0}$ and $u=U^\Gw_{\infty,0}$ we obtain
\bel{uni5}
\lim_{\Gw \ni x\to 0}\myfrac{u^\Gw_{\infty,0}(x)}{U^\Gw_{\infty,0}(x)}=1.
\ee
For any $\vge>0$, by the comparison principle, $(1+\vge)u_{\infty,0}^\Gw \geq U_{\infty,0}^\Gw$ in $\Gw \sms B_\vge$. Letting $\vge \to 0$ yields $u_{\infty,0}^\Gw \geq U_{\infty,0}^\Gw$ in $\Gw$ and thus $u_{\infty,0}^\Gw = U_{\infty,0}^\Gw$ in $\Gw$. \medskip

\noindent{\bf Case 2:} $H$ satisfies \eqref{add} with $p=\frac{q}{2-q}$. The desired result is obtained by a similar argument. \medskip

\noindent{\bf Case 3:} $H$ satisfies \eqref{add} with $p>\frac{q}{2-q}$.
In this case, we use the transformation $t=-\ln r$ for $t\geq 0$ and $ \tilde u(r,\gs)=r^{-\gb_2}v(t,\gs)$ and obtain the following equation instead of \eqref{uni2}
\bel{uni2-add}\BA{lll}
\left(1+\ge^3_1\right)v_{tt}+\left(2\gb_2+2-N+\ge^3_2\right)v_{t}
+\left(\gb_2(\gb_2+2-N)+\ge^3_3\right)v+\Gd'v\\[3mm]
+\langle\nabla'v,\overrightarrow {\ge^3_4}\rangle+\langle\nabla'v_t,\overrightarrow {\ge^3_5}\rangle+
\langle\nabla'\langle \nabla' v,{\bf e}_N\rangle,\overrightarrow {\ge^3_6}\rangle -v^p \\[3mm]
- e^{-\frac{p(2-q)-q}{p-1}t}\abs{(-\gb_1\,v+v_t){\bf n}+\nabla' v+\langle(-\gb_1\,v+v_t){\bf n}+\nabla' v,{\bf e}_N\rangle\overrightarrow\ge^3_7}^q=0
\EA\ee
where $\ge^3_j$ has the same properties as $\ge^1_j$ ($j=\ovl{1,7}$). Notice that $$\lim_{t \to \infty}e^{-\frac{p(2-q)-q}{p-1}t}=0$$ since $p>\frac{q}{2-q}$. By proceeding as in the Case 1, we deduce that $u_{\infty,0}^\Gw = U_{\infty,0}^\Gw$ in $\Gw$. \medskip

\noindent{\bf Case 4:} $H$ satisfies \eqref{add} with $p<\frac{q}{2-q}$. Using a similar argument as in Case 3, we obtain $u_{\infty,0}^\Gw = U_{\infty,0}^\Gw$ in $\Gw$. \qed \medskip

\noindent{\bf Proof of Theorem E.} Statement (i) follows from \rth{UNI}, while statement (ii) follows from \rprop{sing}. \qed
\section{Dirichlet problem with unbounded measure data}
Throughout this subsection we assume that $H$ satisfies \eqref{add}.
\bprop{bound} Assume $H$ satisfies \eqref{add} with $p > 1$, $1< q < 2$ and $K$ is a compact subset of $\prt \Gw$. Then there exists $C>0$ depending on $N$, $p$, $q$ and the $C^2$ characteristic of $\Gw$ such that for any positive solution $u \in C(\ovl \Gw \sms K) \cap C^2(\Gw)$ of \eqref{A0} vanishing on $\prt \Gw \sms K$, there holds
\bel{estK} u(x) \leq C\rho(x)\rho_K(x)^{-\gb_2-1} \forevery x \in \Gw \ee
where $\rho_K(x)=\dist(x,K)$.
\es
\Proof Since $u$ is a positive solution of \eqref{A0}, it is a subsolution of
$$ -\Gd v + v^p = 0 $$
in $\Gw$. By \cite[Proposition 3.4.4]{MVbook}, there exists a constant $C_1>0$ depending on $N$, $p$ and the $C^2$ characteristic of $\Gw$ such that
\bel{ub1} u(x) \leq C_1 \rho(x)\rho_K(x)^{-\frac{p+1}{p-1}} \forevery x \in \Gw. \ee

Next put $\CC_K=\{ x \in \Gw: \rho(x) > \frac{1}{4}\rho_K(x) \}$. Since $u$ is a positive subsolution of
$$ - \Gd v + |\nabla v|^q =0 $$
in $\Gw$, by a similar argument as in the proof of \cite[Proposition 3.5]{NV}, we can show that there exist positive constants $\gd^* \in (0,\gd_0)$ and $C_2>0$ depending on $N$, $q$ and $\Gw$ such that
\bel{ub2} u(x) \leq C_2\rho(x) \rho_K(x)^{-\frac{1}{q-1}} \ee
for every $ x \in \Gw_{\gd^*} \sms \CC_K$. This, along with \eqref{est2-add}, implies that \eqref{ub2} holds in $\Gw$. By combining \eqref{ub1} and \eqref{ub2}, we deduce \eqref{estK}. \qed

\blemma{maximum} Let $u$ and $v$ be two positive solutions of \eqref{A0}. Assume that $u \geq v$ in $\Gw$. Then either $u \equiv v$ or $u>v$ in $\Gw$.
\es
\Proof Put $w=u-v$ then $w \geq 0$ in $\Gw$ and $w$ satisfies
$$-\Gd w + a(x).\nabla w + b(x) w = 0 \q \text{in } \Gw $$
where
$$ a(x)=\left\{ \BA{lll} \frac{(|\nabla u|^q-|\nabla v|^p)\nabla(u-v) }{|\nabla(u-v)|^2} \q &\text{ if } \nabla u \neq \nabla v \\
0 &\text{ if } \nabla u=\nabla v,
\EA \right.$$

$$ b(x)=\left\{ \BA{lll} \frac{u^p-v^p}{u-v} \q &\text{ if } u \neq v \\
0 &\text{ if } u=v
\EA \right.$$
By \rprop{est-add}, $|a(x)| \leq c_1\rho(x)^{-1}$ and $b(x) \leq c_2\rho(x)^{-2}$ in $\Gw$ where $c_1$ and $c_2$ depend on $N$, $p$, $q$ and $\gd_0$.

Next suppose that there exists $x_0 \in \Gw$ such that $w(x_0)=0$. Let $r>0$ such that $B_{3r}(x_0) \sbs \Gw$.
By Harnack inequality \cite[Theorem 5]{Se}, there exists $c_3=c_3(N,p,q,\gd_0,x_0,r)$ such that
$$ \max_{B_\gd(x_0)}w \leq c_3 \min_{B_\gd(x_0)}w = 0. $$
Hence $w \equiv 0$ in $B_r(x_0)$. By standard connectedness argument and Harnack inequality, we deduce that $w \equiv 0$ in $\Gw$. \qed \medskip

For any $k>0$ and $y \in \prt \Gw$, let $u_{k,y}$ be the unique solution of \eqref{A0} with boundary trace $k\gd_y$ and $u_{\infty,y}$ be the unique solution of \eqref{A0} with strong singularity at $y$.  \medskip

\noindent \textbf{Proof of Theorem G.}
\textit{Step 1: Construction of minimal element of $\CU_K$.} Denote by $\CV_K$ the family of all positive moderate solutions $u$ of \eqref{A0} such that $u=0$ on $\prt \Gw \sms K$. Set $u_K:=\sup \CV_K$.

By \rcor{disjoint}, if $u,v\in \CV_K$ then there exists a solution $\tilde u\in \CV_K$ such that $\max(u,v)<\tilde u$.
This fact and \rlemma{maximum} imply (by the same proof as in \cite[Lemma 3.2.1]{MVbook}) that $u_K$ is the limit of an increasing sequence of solutions in $\CV_K$. \rprop{est-add}
implies that $u_K$ is a solution of the equation and vanishes on $\bdw\sms K$.
Clearly $u_K\geq \sup\{u_{\infty,y}: y \in K\}$. Therefore $\CS(u_K)=K$ and $u \in \CU_K$.

Next we show that $u_K$ is the minimal element of $\CU_K$.

If $w \in \CU_K$ then by \rlemma{singusol}
$$ w \geq \sup\{u_{\infty,y}: y \in K\}=\sup\{u_{k,y}: k>0,\; y \in K\}. $$
By Theorem B and \rcor{disjoint}, $w$ dominates every solution of \eqref{A0} whose boundary trace belongs to $\BBD(K)$ ($=$ set of finite linear combination of Dirac measures supported on $K$). If $u\in \CV_K$ then from \rth{moderate} we obtain  $\tr(u)=\mu \in \GTM^+(\prt \Gw)$ with $\supp \mu \sbs K$. Hence there exists a sequence $\{\mu_m\} \sbs \BBD(K)$ converging weakly to $\mu$. By stability and uniqueness result, the  sequence $\{u_{\mu_n}\}$ converges to $u$ in $L^1(\Gw)$. Since $u_{\mu_n} \leq w$ for every $n$, we deduce that $u \leq w$. Therefore $u_K \leq w$ and $u_K$ is the minimal element of $\CU_K$.   \medskip

\noindent \textit{Step 2: Construction of maximal element of $\CU_K$.}  Denote by $\CW_K$ the family of all positive solutions $u$ of \eqref{A0} such that $u=0$ on $\prt \Gw \sms K$. Put $U_K:=\sup \CW_K$. By the same argument as in Step 1, one shows that $U_K \in \CW_K$. By \rlemma{singusol}, $ U_K \geq \sup\{u_{\infty,y}: y \in K\}$,
which implies $\CS(U_K)=K$. Therefore  $U_K$ is the maximal element of $\CU_K$.  \medskip

\noindent \textit{Step 3: Proof of \eqref{2K}.} Pick $y \in K$. We may assume $y$ is the origin. By \rprop{sing}, for every $\gg \in (0,1)$, there exists $r=r(\gg)$ and $c=c(N,p,q,\gg)$ such that
\bel{j1} u^\Gw_{\infty,y}(x) \geq c |x-y|^{-\gb_2} \forevery x \in C_{\gg,r}(y):=\{x \in \Gw: \rho(x) \geq \gg |x-y|\} \cap B_r(y). \ee
Since $y \in K=\CS(u_K)$ we have $u_K \geq u_{\infty,y}$. Therefore
\bel{j2} u_K(x) \geq c |x-y|^{-\gb_2} \forevery x \in C_{\gg,r}(y). \ee
On the other hand, by \rprop{est-add},  for every $x \in C_{\gg,r}(y)$,
\bel{j3} U_K(x) \leq \tl \Gl_2 \rho(x)^{-\gb_2} \leq \tl \Gl_2 \gg^{-\gb_2}|x-y|^{-\gb_2}. \ee
From \eqref{j2} and \eqref{j3} we deduce \eqref{2K}. \qed

\section{Removability}
In this section we deal with removable singularities in the case that $H$ is supercritical.
\bprop{nostrong} Assume either $H$ satisfies \eqref{multi} with $N(p+q-1)\geq p+1$
or H satisfies \eqref{add} with $m_{p,q} \geq p_c$. If  $u\in C(\overline\Gw\setminus\{0\})\cap C^2(\Gw)$ is a nonnegative solution of \eqref{A0} vanishing on $\prt\Gw\setminus\{0\}$ then $u$ cannot be a strongly singular solution. \es
\Proof  We consider a sequence of functions $\gz_n\in C^\infty(\BBR^N)$ such that $\gz_n(x)=0$ if $|x|\leq \frac{1}{n}$, $\gz_n(x)=1$ if $|x|\geq \frac{2}{n}$, $0\leq \gz_n\leq 1$ and $|\nabla\gz_n|\leq c_{13}n$, $|\Gd\gz_n|\leq c_{13}n^2$ where $c_{13}$ is independent of $n$. We take $\gx\gz_n$ as a test function (where $\gx$ is the solution to \eqref{eta}) and we obtain
\bel{est3}\BA {l}
\myint{\Gw}{}( u + (H \circ u)\gx)\gz_n \, dx=\myint{\Gw}{}u\left(\gx\Gd \gz_n+2\nabla\gx.\nabla\gz_n \right)dx
=:J+J'.
\EA\ee
Set $\CO_n=\Gw\cap \{x:\frac{1}{n}<|x|\leq \frac{2}{n}\}$, then $|\CO_n|\leq c_{14}(N)n^{-N}$. On the one hand, since $\gx(x)\leq c_3\rho(x) \leq c_3|x|$,
$$  J \leq c_{15}\Gl_i\,\myint{\CO_n}{}n^{\gb_i+2}\gx dx\leq c_{16}n^{\gb_i+1-N}
$$
 where
\bel{i} i=\left\{ \BA{lll} 1 \q \text{if } H \text{ satisfies } (\ref{multi}), \\ 2 \q \text{if } H \text{ satisfies } (\ref{add}). \EA \right. \ee
On the other hand,
\bel{est3'} J' \leq c_{17}\Gl_i\,\myint{\CO_n}{}n^{\gb_i+1}|\nabla\gx| dx\leq c_{18}n^{\gb_i+1-N} \ee
where $i$ is given in \eqref{i}. By combining \eqref{est3}-\eqref{est3'} and then by letting $n \to \infty$ we obtain
\bel{est4}\BA {l}
\myint{\Gw}{}\left(u+ (H \circ u)\gx\right) dx < \infty. \EA\ee
By \rth{moderate}, the boundary trace of $u$ is a finite measure. Since $u=0$ on $\prt \Gw \sms \{0\}$, the boundary trace of $u$ is $k\gd_0$ for some $k \geq 0$. \qed

\bcor{supercritical} Assume either $H$ satisfies \eqref{multi} with $N(p+q-1) > p+1$
or H satisfies \eqref{add} with $m_{p,q} > p_c$. If  $u\in C(\overline\Gw\setminus\{0\})\cap C^2(\Gw)$ is a nonnegative solution of \eqref{A0} vanishing on $\prt\Gw\setminus\{0\}$ then $u \equiv 0$. \es
\Proof Since $\gb_i+1-N<0$, we deduce from \eqref{est3}-\eqref{est3'} that
$$ \myint{\Gw}{}(u+(H \circ u)\gx) dx = 0, $$
which implies $u \equiv 0$. \qed

\bth{critical} Assume either $H$ satisfies \eqref{multi} with $N(p+q-1) = p+1$
or H satisfies \eqref{add} with $m_{p,q} = p_c$. If $u\in C(\overline\Gw\setminus\{0\})\cap C^2(\Gw)$ is a nonnegative solution of \eqref{A0} vanishing on $\prt\Gw\setminus\{0\}$ then $u\equiv 0$.
\es
\Proof By \rprop{nostrong}, $u$ admits a boundary trace  $k\gd_0$, $k \geq 0$.

 For $0<\ell<1$, we set
    $$ u_\ell(x)=T^1_\ell[u](x)=T^2_\ell[u](x)=\ell^{N-1}u(\ell x), \qq x\in \Gw^{\ell}=\ell^{-1}\Gw. $$
By the comparison principle, $u_\ell \leq kP^{\Gw^\ell}(.,0)$ in $\Gw^\ell$  for every $\ell \in (0,1)$. Due to \rprop{halfspace-c}, up to a subsequence, $\{u_\ell\}$ converges to a function $\tl u$ which is a solution of either \eqref{half1} if $H$ satisfies \eqref{multi}, or \eqref{half2} if $H$ satisfies \eqref{add} with $p=\frac{q}{2-q}$, or \eqref{half3} if $H$ satisfies \eqref{add} with $p>\frac{q}{2-q}$, or \eqref{half4} if $H$ satisfies \eqref{add} with $p<\frac{q}{2-q}$. Moreover, $\tl u \leq kP^{\BBR^N_+}(.,0)$ in $\BBR^N_+$. \medskip

If $H$ satisfies \eqref{add} with $p \neq \frac{q}{2-q}$ then since $m_{p,q}=p_c$, it follows from \cite{MV1} and \cite{NV} that $\tl u=0$. \medskip

If $H$ satisfies \eqref{multi} or $H$ satisfies \eqref{add} with $p = \frac{q}{2-q}$ then set
$$\CV=\{ v: v \text{ is a solution of } \eqref{A0} \text{ in } \BBR^N_+, \tl u \leq v \leq kP^{\BBR^N_+}(.,0) \}$$
and put $\tl v:=\sup\CV$. \medskip

\noindent{\it Assertion: } $\tl v$ is a solution of \eqref{halfspace} in $\BBR^N_+$.

Indeed, let $\{Q_n\}$ be a sequence of $C^2$ bounded domains such that $\ovl Q_n \sbs Q_{n+1}$, $\cup_{n \in \BBN}Q_n=\BBR^N_+$ and $0<\dist(Q_n,\prt \BBR^N_+)<\frac{1}{n}$ for each $n \in \BBN$. Consider the problem
    \bel{Dn} \left\{ \BA{lll} -\Gd w +  H \circ w & = 0 \qq &\text{in } Q_n \\
    \phantom{-\Gd w +  H \circ ,}
    w&= kP^{\BBR^N_+}(.,0) &\text{on } \prt Q_n.\EA \right. \ee
Since $\tl u$ and $kP^{\BBR^N_+}(.,0)$ are respectively subsolution and supersolution of \eqref{Dn}, there exists a solution $w_n$ of the problem \eqref{Dn} satisfying $\tl u \leq w_n \leq kP^{\BBR^N_+}(.,0)$ in $Q_n$. Hence, by the comparison principle $\tl u \leq w_{n+1} \leq w_n \leq kP^{\BBR^N_+}(.,0)$ in $Q_n$ for each $n \in \BBN$. Therefore, $\tl w:=\lim_{n \to \infty}w_n \leq kP^{\BBR^N_+}(.,0)$ in $\BBR^N_+$. By regularity results \cite{Lib}, we obtain \eqref{lie} with $v_\ell$ replaced by $w_n$ and $\Gw^\ell$ replaced by $Q_n$. Thus $\tl w$ is a solution of \eqref{halfspace}. On the one hand, by the definition of $\tl v$, $\tl w \leq \tl v$. On the other hand, $\tl v \leq w_n$ in $Q_n$ for every $n$, and consequently $\tl v \leq \tl w$ in $\BBR^N_+$. Thus $\tl v= \tl w$. \medskip

For every $\ell>0$, we set $w_\ell=T^1_\ell[\tl v]=T^2_\ell[\tl v]=\ell^{N-1}\tl v(\ell x)$ with $x \in \BBR^N_+$ then $w_\ell=\sup\CV=\tl v$ in $\BBR^N_+$ for every $\ell>0$. Hence $\tl v$ is self-similar, namely $\tl v$ can be written under the separable form
    $$ \tl v(r,\gs)=r^{N-1}\gw_i(\gs) \qq (r,\gs) \in (0,\infty) \ti S^{N-1}_+$$
where $\gw_i$ is the nonnegative solution of \eqref{PH1}. It follows from \rth{exis-uniq-SS} that $\gw_i \equiv 0$, hence $\tl v \equiv 0$. Thus $\tl u \equiv 0$.

Hence
    \bel{t4} \lim_{n \to \infty}(\sup\{\abs{ u_{\ell_n}(x)}+\abs{\nabla  u_{\ell_n}(x)}: x \in \Gg_{R^{-1},R} \cap \Gw^{\ell_n}\})=0. \ee
Consequently,
    $$ \lim_{x \to 0}\abs{x}^{N-1}u(x)=0 \qq \text{and} \qq \lim_{x \to 0}\abs{x}^N\abs{\nabla u(x)}=0. $$
Therefore, $ \lim_{x \to 0}(\abs{x}^N\rho(x)^{-1}u(x))=0$, namely $u=o(P^\Gw(.,0))$. By the comparison principle, $u \equiv 0$. \qed    \medskip

\noindent{\bf Proof of Theorem F.} The proof follows immediately from \rcor{supercritical} and \rth{critical}. \qed

We next deal with the case $q=2$.

\bth{q=2} Assume $q=2$.  If $u\in C(\overline\Gw\setminus\{0\})\cap C^2(\Gw)$ is a nonnegative solution of \eqref{A0}
 vanishing on $\prt\Gw\setminus\{0\}$ then $u\equiv 0$. \es
\Proof Put
$$v=\left\{ \BA{lll}1-e^{-\frac{1}{p+1}u^{p+1}} &\text{ if } H \text{ satisfies } (\ref{multi}), \\
1-e^{-u} &\text{ if } H \text{ satisfies } (\ref{add}) \EA \right. $$
then $v \in C(\overline\Gw\setminus\{0\})\cap C^2(\Gw)$, $0 \leq v \leq 1$ and $v$ satisfies
	\bel{v} -\Gd v \leq 0 \q \text{in } \Gw, \qq v=0 \q \text{on } \prt \Gw \sms \{0\}. \ee
Let $\eta_\gd$ be the solution of
	\bel{vd} -\Gd \eta_\gd = 0 \q \text{in } D_\gd, \qq \eta_\gd=v \q \text{on } \prt D_\gd \ee
then by the comparison principle  $v \leq \eta_\gd \leq 1$ in $D_\gd$. The sequence $\{\eta_\gd\}$ converges to an harmonic function $\eta^* \geq v$ as $\gd \to 0$. Since $0 \leq \eta^* \leq 1$ and $\eta^*=0$ on $\prt \Gw \sms \{0\}$, it follows that $\eta^* \equiv 0$. Hence $v \equiv 0$, so $u \equiv 0$. \qed \medskip

\appendix \section{Uniqueness result in subcritical case
   \\ by Phuoc-Tai Nguyen}
\setcounter{equation}{0}

In this section, we deal with the question of uniqueness for the problem \eqref{A1}. Let $\Gw$ be a $C^2$ bounded domain in $\BBR^N$. We assume that $H \in C(\Gw \times \BBR \times \BBR^N)$ satisfies
\bel{asH} |H(x,u,\xi)-H(x,u',\xi')| \leq A\rho(x)^{\ga}(a(x)+|\xi|^{q-1}+|\xi'|^{q-1})|\xi -\xi'| \ee
for a.e. $x \in \Gw$, for every $u,u' \in \BBR$ and $\xi, \xi' \in \BBR^N$, where  $A>0$, $\ga \in (-1,\frac{1}{N-1})$, $q \in (1, q_{\ga,c})$ with $q_{\ga,c}:=\frac{N+1+\ga}{N}$  and  $a \in L_{\rho^{1+\ga}}^{\frac{q}{q-1}}(\Gw)$. As above, we use the notation $H\circ u$ to denote $H(x,u(x),\!\nabla u(x))$.

Solutions of \eqref{A1} are always understood in the sense of \rdef{sol}.

The uniqueness result is stated as follows:
\bth{uniqueA} Assume $H$ satisfies \eqref{asH}. For every $\gm \in \GTM^+(\prt \Gw)$, \eqref{A1} admits at most one solution. \es
The proof of \rth{uniqueA} is an adaptation of the method in \cite{Po} and based upon the following lemma
\blemma{lem2} Let $f \in L^1_\rho(\Gw)$ and $z$ be a positive solution of
\bel{b8} - \Gd z \leq f \quad \text{in } \Gw, \quad z=0 \quad \text{on } \prt \Gw. \ee
Then for any $\gg \in (0,\frac{N}{N-1})$ and $1<q<\frac{N+\gg}{N}$, there exists a constant  $c_{19}=c_{19}(N,\Gw,\gg)$ such that
$$\norm{\nabla z}_{L^q_{\rho^\gg}(\Gw)} \leq c\norm{f}_{L^1_\rho(\Gw)}. $$
\es
\Proof Let $\tl z$ is the unique solution to the problem
\bel{b8'} - \Gd \tl z = f \quad \text{in } \Gw, \quad \tl z=0 \quad \text{on } \prt \Gw. \ee
then there exists a positive constant $c_{20}=c_{20}(N,\Gw,\gg)$ such that
	\bel{b8*} \norm{\nabla \tl z}_{L^q_{\rho^\gg}(\Gw)} \leq c_{20}\norm{f}_{L^1_\rho(\Gw)}.\ee
Denote $\eta = \tl z -z$ then $\eta$ satisfies
	$$ - \Gd \eta \geq 0 \quad \text{in } \Gw, \quad \eta=0 \quad \text{on } \prt \Gw. $$
By the maximum principle, $\eta \geq 0$ in $\Gw$. As $\eta$ is a superharmonic function, by \cite[Theorem 1.4.1]{MVbook}, there exists a positive Radon measure $\tau \in \GTM_{\rho}^+(\Gw)$ such that
	\bel{b8a} - \Gd \eta = \tau \quad \text{in } \Gw, \quad \eta=0 \quad \text{on } \prt \Gw. \ee
 Therefore, by \cite{BVi}, for every $q \in (1,\frac{N+\gg}{N})$, $|\nabla \eta| \in L^q_{\rho^\gg}(\Gw)$  and
	\bel{b8b}\norm{\nabla \eta}_{L^q_{\rho^\gg}(\Gw)} \leq c_{21}\norm{\tau}_{\GTM_\rho(\Gw)} \ee
where $c_{21}$ is another positive constant depending on $N$, $\Gw$ and $\gg$. Moreover, by using $\xi$ ($\xi$ is the solution of \eqref{eta}) as a test function, we deduce from \eqref{b8a} that
	$$ \norm{\eta}_{L^1(\Gw)} = \int_{\Gw}\xi d \tau. $$
From \eqref{b8*} and \eqref{b8b}, $z=\tl z -\eta \in L^q_{\rho^\gg}(\Gw)$ and
	\bel{b8c}  \norm{\nabla z}_{L^q_{\rho^\gg}(\Gw)} \leq c_{22} (\norm{f}_{L^1_\rho(\Gw)} + \norm{\tau}_{\GTM_\rho(\Gw)}) \ee
where $c_{22}=c_{22}(N,\gg,\Gw)$. Since $z \geq 0$, it follows that
	\bel{b8d} \norm{\tau}_{\GTM_\rho(\Gw)} \leq c_3  \int_{\Gw}\xi d \tau = c_3\norm{\eta}_{L^1(\Gw)} \leq c\norm{\tl z}_{L^1(\Gw)} \leq c_{23}\norm{f}_{L^1_\rho(\Gw)} \ee
where $c_{23}=c_{23}(N,\gg,\Gw)$. From \eqref{b8c} and \eqref{b8d} we get the desired estimate. \qed \medskip

We turn to the

\noindent{\bf Proof of \rth{uniqueA}.} Let $u_1$ and $u_2$ be two solutions of \eqref{A1}. Since $q<q_{\ga,c}$, from \rprop{P1}, we deduce that $|\nabla u_i| \in L^q_{\rho^{1+\ga}}(\Gw)$, $i=1,2$  and
$$ \norm{\nabla u_i}_{L^q_{\rho^{1+\ga}}(\Gw)} \leq c_1(\norm{H \circ u_i}_{L^1_{\rho}(\Gw)}+ \norm{\gm}_{\GTM(\prt \Gw)}).
$$
Let $\{\gm_n\}$ be a sequence of functions in $C^1(\prt \Gw)$ converging weakly to $\gm$. For $k>0$, denote by $T_k$ the truncation function, i.e. $T_k(s)=\max(-k,\min(s,k))$. For very $n>0$, denote by $u_{i,n}$, $i=1,2$ the solution of the problem
\bel{b2} - \Gd u_{i,n} + T_n(H \circ u_i) = 0 \quad \text{in } \Gw, \quad u_{i,n}=\gm_n \quad \text{on } \prt \Gw.
\ee
By local regularity theory for elliptic equations (see, e.g., \cite{Mi}), $u_{i,n} \to u_i$ in $C^1_{loc}(\Gw)$ for $i=1,2$.
From \eqref{b2} we obtain
\bel{b6} \left\{ \BA{lll} - \Gd (u_{1,n}-u_{2,n}) &= -T_n(H \circ u_1) + T_n(H \circ u_2) \quad \text{in } \Gw \\
\phantom{--,,}
 u_{1,n}-u_{2,n} &= 0 \quad \text{on } \prt \Gw. \EA \right. \ee
We shall prove that $u_1 \leq u_2$. By contradiction, we assume that $M:=\sup_\Gw(u_1-u_2) \in (0,\infty]$. Let $0<k<M$. From \eqref{b6} and Kato's inequality \cite{MVbook}, we get
\bel{b7} -\Delta(u_{1,n}-u_{2,n}-k)_+ \leq (T_n(H \circ u_2) - T_n(H \circ u_1))\chi_{_{E_{n,k}}} \ee
where $E_{n,k}=\{x \in \Gw: u_{1,n}-u_{2,n}>k\}$. Applying \rlemma{lem2} with $\gg=1+\ga$ and Holder's inequality, thanks to \eqref{asH}, we get
\bel{b9} \BA{lll} \myint{\Gw}{}|\nabla(u_{1,n}-u_{2,n}-k)_+ |^q\rho^{1+\ga} dx \\
\phantom{}
\leq c_{24}\left(\myint{\Gw}{}|(T_n(H \circ u_2) - T_n(H \circ u_1)|\chi_{_{E_{n,k}}}\rho dx\right)^q\\
\phantom{}
\leq c_{24}A^q\left(\myint{F_{n,k}}{}(a(x)+|\nabla u_1|^{q-1}+|\nabla u_2|^{q-1})|\nabla(u_1-u_2)|\rho^{1+\ga} dx\right)^q\\
\phantom{}
\leq c_{25}\left(\myint{F_{n,k}}{}(a(x)^{\frac{q}{q-1}} +|\nabla u_1|^q+ |\nabla u_2|^q)\rho^{1+\ga} dx\right)^{q-1}\myint{F_{n,k}}{}|\nabla(u_1-u_2)|^q\rho^{1+\ga} dx
\EA \ee
where $F_{n,k}=\{x \in \Gw: u_{1,n}-u_{2,n}>k, \nabla u_1 \neq \nabla u_2\}$. Since $u_{1,n}-u_{2,n} \to u_1-u_2$ a.e. in $\Gw$, $\chi_{_{F_{n,k}}} \to \chi_{_{F_k}}$ a.e. where $F_k=\{x \in \Gw: u_1-u_2>k, \nabla u_1 \neq \nabla u_2\}$. Hence, letting $n \to \infty$ in \eqref{b9} implies
\bel{b10} \myint{\Gw}{}|\nabla(u_1-u_2-k)_+ |^q\rho^{1+\ga} dx \leq c_{25}R_k \myint{F_{k}}{}|\nabla(u_1-u_2-k)_+|^q\rho^{1+\ga} dx \ee
where
$$ R_k:= \left(\myint{F_{k}}{}(a(x)^{\frac{q}{q-1}}+|\nabla u_1|^q+ |\nabla u_2|^q)\rho^{1+\ga} dx\right)^{q-1}. $$
\noindent{\bf Case 1: $M=\infty$.} Then $\lim_{k \to \infty}|F_k|=0$. Since $a \in L_{\rho^{1+\ga}}^{\frac{q}{q-1}}(\Gw)$ and $|\nabla u_1|, |\nabla u_2| \in L^q_{\rho^{1+\ga}}(\Gw)$, there exists $k_0$ large enough so that $R_{k_0} < (2c_{25})^{-1}$. Hence, \eqref{b10} implies
\bel{b11} \myint{\Gw}{}|\nabla(u_1-u_2-k_0)_+ |^q\rho^{1+\ga} dx=0. \ee
It follows that $\nabla(u_1-u_2-k_0)_+=0$ in $\Gw$ and hence $(u_1-u_2-k_0)_+=c_0 \geq 0$. Therefore, $u_1-u_2 \leq k_0+c_0$ a.e. in $\Gw$, which contradicts $\sup_\Gw(u_1-u_2)=\infty$. \medskip

\noindent{\bf Case 2: $M<\infty$}. Since $|\nabla (u_1-u_2)|=0$ a.e. on the set $\{x \in \Gw: (u_1-u_2)(x)=M  \}$, it follows that $\lim_{k \to M}|F_k|=0$. By proceeding as in Step 1, we deduce that there exists $k_0 \in (0,M)$ and $c_0 \geq 0$ such that $(u_1-u_2-k_0)_+ = c_0$ a.e. in $\Gw$. If $c_0=0$ then $u_1-u_2 \leq k_0$, which contradicts $\sup_\Gw(u_1-u_2)=M>k_0$. If $c_0>0$ then $u_1-u_2=k_0+c_0$. Then $u_1-u_2=M$ a.e. in $\Gw$, which contradicts
the fact that $u$ and $u_2$ have the same boundary trace $\gm$.

Thus $u_1 \leq u_2$ and the uniqueness follows by permuting the roles of $u_1$ and $u_2$.\qed \bigskip

\noindent{\bf Acknowledgements} The authors were supported by the Israel Science Foundation founded by the Israel Academy of Sciences and Humanities, through grant 91/10. The second author was supported in part by a Technion fellowship. The second author wishes to thank Prof. L. V\'eron and Prof. A. Porretta for fruitful discussion concerning results in Appendix A.

\end{document}